\documentclass{amsart}
\usepackage{amssymb}
\usepackage{mathtools}
\usepackage{mathrsfs}
\usepackage[hidelinks]{hyperref}
\usepackage{tikz}
\usepackage{graphicx}
\usepackage{todonotes}
\usepackage[font=small]{caption}
\usetikzlibrary{decorations.markings}

\tikzset{
	wing/.style={
		decoration={
			markings, mark=at position 0.8 with {\arrow{latex}}
		},
		postaction={decorate},
	},
}

\newcommand{\Z}{\mathbb{Z}}
\newcommand{\R}{\mathbb{R}}
\newcommand{\C}{\mathbb{C}}

\newcommand{\ii}{{\rm i}}
\newcommand{\cent}[1]{#1^\circ}

\newcommand{\graph}{\mathsf}
\newcommand{\gV}{\graph{V}}
\newcommand{\gE}{\graph{E}}
\newcommand{\gR}{\graph{R}}
\newcommand{\gF}{\graph{F}}
\newcommand{\gH}{\graph{H}}
\newcommand{\gB}{\graph{B}}
\newcommand{\gC}{\graph{C}}

\newcommand{\surface}{\mathscr}
\newcommand{\sM}{\surface{M}}
\newcommand{\sS}{\surface{S}}

\newcommand{\Fhor}{F^\text{hor}}
\newcommand{\Phor}{P^\text{hor}}
\newcommand{\Fver}{F^\text{ver}}
\newcommand{\Pver}{P^\text{ver}}

\newcommand{\cA}{\mathcal{A}}
\newcommand{\cB}{\mathcal{B}}
\newcommand{\cC}{\mathcal{C}}

\newcommand{\cR}{\mathcal{R}}
\newcommand{\cD}{\mathcal{D}}
\newcommand{\cM}{\mathcal{M}}

\newcommand{\cPhor}{\mathcal{P}^\text{hor}}
\newcommand{\cFhor}{\mathcal{F}^\text{hor}}
\newcommand{\cPver}{\mathcal{P}^\text{ver}}
\newcommand{\cFver}{\mathcal{F}^\text{ver}}

\newcommand{\rotate}{\varsigma}
\newcommand{\involute}{\iota}

\renewcommand{\div}{\operatorname{div}}
\newcommand{\curl}{\operatorname{curl}}
\newcommand{\mdiv}{\operatorname{mdiv}}

\newcommand{\re}{\operatorname{Re}}

\newcommand{\Res}[2]{\operatorname{Res}\left(#1,#2\right)}

\theoremstyle{plain} 
\newtheorem{theorem}{Theorem}
\newtheorem{proposition}{Proposition}
\newtheorem{lemma}[proposition]{Lemma}

\theoremstyle{definition} 
\newtheorem{definition}{Definition}
\newtheorem{example}{Example}

\theoremstyle{remark}
\newtheorem{remark}{Remark}

\begin{document}

\title[Gluing saddle towers II: SPMS]{Gluing Karcher--Scherk saddle towers II:\\Singly periodic minimal surfaces}

\author{Hao Chen}
\address[Chen]{ShanghaiTech University, Institute of Mathematical Sciences}
\email{chenhao5@shanghaitech.edu.cn}
\thanks{H.\ Chen was partially supported by Individual Research Grant from Deutsche Forschungsgemeinschaft within the project ``Defects in Triply Periodic Minimal Surfaces'', Projektnummer 398759432.}


\date{\today}

\begin{abstract}
	This is the second in a series of papers that construct minimal surfaces by
	gluing singly periodic Karcher--Scherk saddle towers along their wings.  This
	paper aims to construct singly periodic minimal surfaces with Scherk ends.
	As in the first paper, we prescribe phase differences between saddle towers,
	and obtain many new examples without any horizontal reflection plane.  This
	construction is not very different from previous ones, hence we will only
	provide sketched proofs.

	We will however study the embeddedness with great care.  Previously,
	embeddedness can not be determined in the presence of ``parallel'' Scherk
	ends, as it was not clear if they bend towards or away from each other.  In a
	recent study, the bending was completely ignored and embeddedness was falsely
	claimed.  We correct this mistake by carefully analysing slight bendings,
	thus identify scenarios where the constructed surfaces are indeed embedded.
\end{abstract}

\maketitle

\section{Introduction}

In~\cite{traizet2001}, Traizet desingularized arrangements of vertical planes
into singly periodic minimal surfaces (SPMSs) using the node-opening
technique\footnotemark: Scherk towers are placed at the intersection lines and
are glued along their wings.  However, the construction relied on many
assumptions: (1) The arrangement was assumed to be \emph{simple} in the sense
that no three planes intersect in a line; (2) The minimal surface was assumed
to be symmetric in a horizontal plane; (3) For the surfaces to be embedded, it
was assumed that no two planes are parallel.  The purpose of this paper is to
get rid of these assumptions as far as we can.

\footnotetext{The construction was also implemented earlier by solving
non-linear PDEs~\cite{traizet1996}.}

On the one hand, we will consider a larger family of configurations.  More
specifically, let $\Gamma$ be a ``graph'' which, informally speaking, consists
of straight segments and rays (edges) that intersect only at their endpoints
(vertices).  We will desingularize $\Gamma \times \R$ into a minimal surface by
placing Karcher--Scherk saddle towers over the vertices and glue them along
their wings following the pattern of the graph.  In particular, the minimal
surface has Scherk ends corresponding to the rays of the graph.

On the other hand, we will remove the horizontal reflection plane.  This is
done, as in the first paper, by prescribing phase differences between saddle
towers.  Our main result (Theorem~\ref{thm:main} below) is then analogous to
that of the first paper~\cite{saddle1}: \emph{The gluing construction sketched
above produces a continuous family of immersed SPMSs only when the graph
satisfies a horizontal balancing condition, and the phases of the saddle towers
satisfy a subtle vertical balancing condition}.  Consequently, we obtain many
examples without any reflectional symmetry; see Section~\ref{ssec:gyroid}.  They
provide a negative answer to a question in~\cite{traizet1996} that asks whether
every SPMS with Scherk ends has a horizontal reflection plane.

\begin{remark}
	The first SPMSs with Scherk ends but no horizontal reflection plane were
	explicitly constructed in~\cite{martin2006} with a very different technique.
	Our construction demonstrates that Traizet's node-opening
	technique~\cite{traizet2001} is flexible and powerful in producing implicit
	non-symmetric examples.
\end{remark}

In this paper, we will only sketch the gluing construction, as all technical
details can be found in the first paper~\cite{saddle1} or even earlier works,
and we do not want to repeat ourselves.  The readers are therefore expected to
have a reasonable familiarity with the first paper.

\medskip

Our main concern is the embeddedness, especially when the graph $\Gamma$ has
parallel rays.  In~\cite{traizet2001}, the embeddedness was only guaranteed in
the absence of parallel vertical planes because, otherwise, the corresponding
Scherk ends risk to bend towards each other after desingularization, therefore
create self-intersection.

In some recent work~\cite{morabito2020}, the bendings of Scherk ends were
completely ignored and embedded SPMSs with ``parallel'' Scherk ends were
\emph{falsely} claimed.  We are therefore compelled to provide a proper
technical treatment on the bendings.  Indeed, our construction allows
quantitative detections of very slight bendings, thus helps to resolve very
delicate embeddedness.

\begin{figure}[!hbt]
	\includegraphics[width=.9\textwidth]{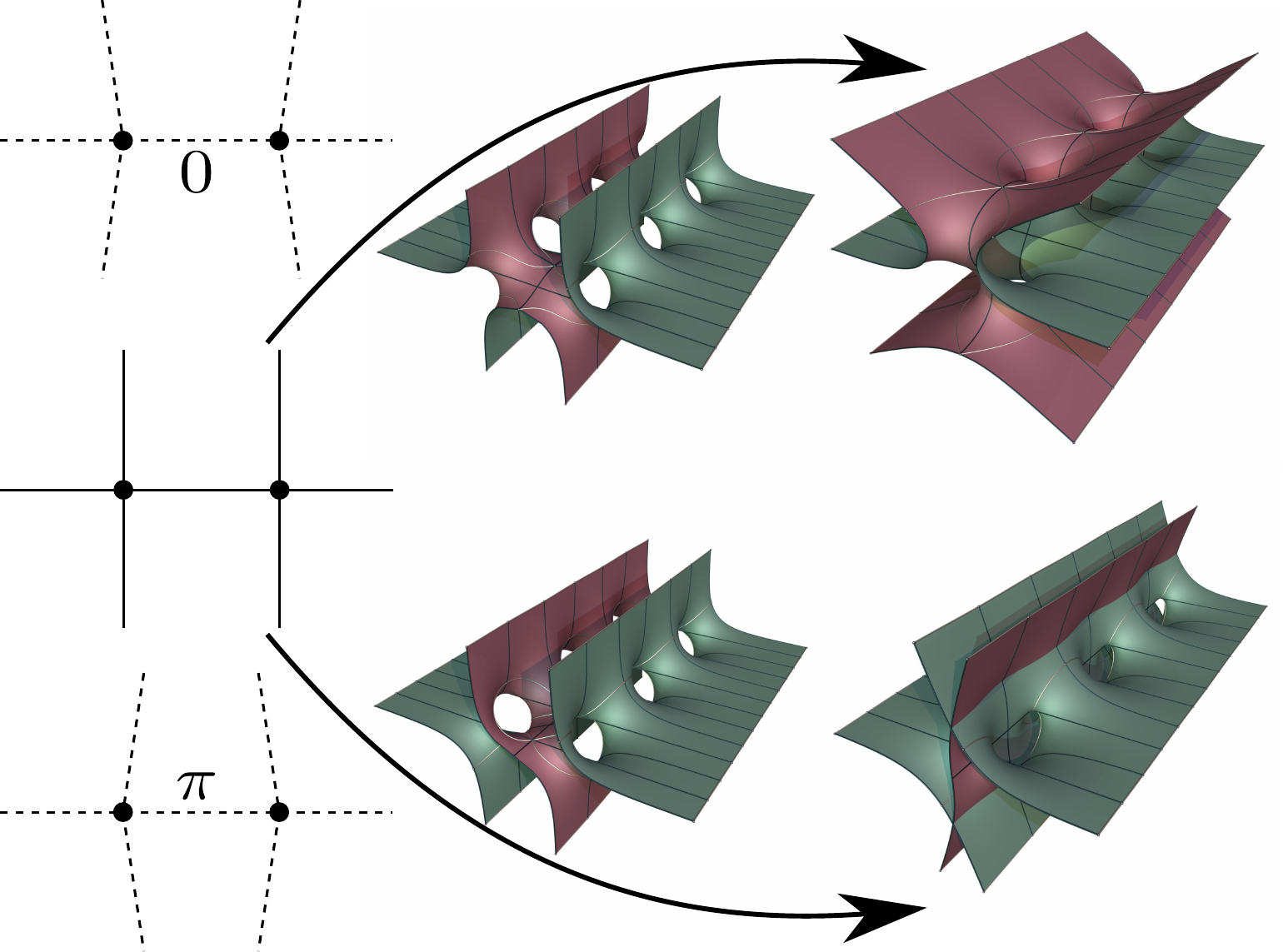}
	\caption{
		The simple graph on the left gives rise to SPMSs with two saddle towers
		that are glued either in phase or in opposite phases.  If they are in
		phase, the surfaces are embedded after desingularization, as shown in the
		figures on the top.  If they are in opposite phases, the surfaces are not
		embedded, as shown in the figures in the bottom.
		\label{fig:tree1}
	}
\end{figure}

We identify two types of bending.  The first arise from the horizontal
deformations, as the saddle towers expand and the glued wings shrink.  We will
explicitly describe this deformation to the lowest order; see
Theorem~\ref{thm:embed2}.  However, in the case of simple vertical plane
arrangment~\cite{traizet2001}, for example, this deformation does not help
determine embeddedness.

We are then obliged to consider a very subtle type of bending, arising from the
need to balance slight variations of the horizontal forces.  This bending is
very delicate.  We will see that, while the expansion of saddle towers and
shrinking of the glued wings are of the order $\varepsilon^2$, the variations
of the horizontal forces, as well as the deformations they cause, are of the
order $\exp(-\ell/\varepsilon^2)$, where $\ell$ is the shortest edge length in
the graph.  It is therefore understandable that these bendings could be easily
ignored.  Again we will explicitly describe this deformation to the lowest
order; see Theorem~\ref{thm:embed4}.

\begin{example}\label{ex:3line}
	Figure~\ref{fig:tree1} shows the simplest example that demonstrates the
	subtle bending.  The graph appears as an arrangement of three lines, one
	horizontal and two vertical.  So we will glue two Scherk saddle towers along
	a single pair of wings.  We will see in Lemma~\ref{lem:tree} that the two
	saddle towers are either in phase or in opposite phases.  Contrary to what
	was claimed in~\cite{morabito2020}, the Scherk ends will not remain parallel
	after desingularization.  If the saddle towers are in phase, the force along
	the glued wings will increase, and the Scherk ends must bend away from each
	other to preserve balance, so the resulting surfaces are embedded.  If they
	are in opposite phases, the force will decrease, and the Scherk ends must
	bend towards each other, so the resulting surfaces are not embedded.  We will
	revisit this example in Section~\ref{ssec:tree}. \qed
\end{example}

\medskip

The paper is organized as follows.  In Section~\ref{sec:result}, we set up the
graph theoretical language before using it to state our main results.  Examples
are given in Section~\ref{sec:examples}.  Section~\ref{sec:proof} is dedicated
to the constructions and proofs.

The construction of the immersed families will only be sketched, as the
technical details can be found in the first paper of this
series~\cite{saddle1}.  Only the proof of embeddedness
(Section~\ref{ssec:embed}), especially in the case of simple vertical plane
arrangement (Section~\ref{sss:subtle}), will be given in detail, because the
involved technique can not be found in previous works.

\subsection*{Acknowledgement}

I appreciate the quick and friendly response from Fillippo Morabito upon
learning Example~\ref{ex:3line}.   He has acknowledged the existence of mistake
in~\cite{morabito2020}.

All 3D pictures in this paper are from \url{http://minimalsurfaces.blog}, an
online repository maintained by Matthias Weber, to whom I express my gratitude.
I also thank Peter Connor who pointed me to some known examples that arise from
our construction.

\section{Main result}\label{sec:result}

\subsection{Graph theory}

We define a \emph{pseudo rotation system} as a triplet $(\gH, \involute,
\rotate)$ where $\gH$ is a finite set of \emph{half-edges}, $\involute$ and
$\rotate$ are two permutations acting on $\gH$, $\involute$ is an involution,
and the group generated by $\involute$ and $\rotate$ acts transitively on
$\gH$.  Note that $(\gH, \involute, \rotate)$ is not a rotation system (hence
``pseudo'') because we allow the involution $\involute$ to have fixed points.
We use $\gR$ to denote the set of fixed points of $\involute$.

In analogy to the rotation systems that define multigraphs, a pseudo rotation
system defines a graph-like structure $(\gH, \gV, \gE)$, where the vertex set
$\gV$ consists of the orbits of $\rotate$, and the edge set $\gE$ consists of
the orbits of $\involute$.  Edges with two half-edges are called \emph{closed
edges}; they are like the edges in the traditional sense.  But we also have
edges with single half-edge; they are called \emph{open edges}, and are
identified with the fixed points of $\involute$.

\begin{remark}[Notation]
	In the remaining of the paper, we will use the letters $h$, $r$, and $\eta$
	to denote the half-edges in, respectively, $\gH \setminus \gR$, $\gR$, and
	$\gH$.  For each $\eta \in \gH$, we use $v(\eta)$ and $e(\eta)$ to denote the
	unique vertex and edge associated to $\eta$.  For a half-edge $h \in \gH
	\setminus \gR$, we write $-h$ for $\iota(h)$; this notation does not apply to
	half-edges $r \in \gR$.  The cardinality of a vertex $v$, seen as a set of
	half-edges, is the degree of the vertex, denoted by $\deg(v)$.
\end{remark}

We assume that the structure admits a \emph{geometric representation} $\varrho$
that maps vertices to distinct points in $\R^2 \simeq \C$, closed edges to line
segments, and open edges to rays, so that the image of each edge is bounded by
the images of its end vertices, and the image of different edges are either
identical or have disjoint interiors.  Closed edges with the same image are
called \emph{parallel}; parallelism is an equivalence relation on the set of
closed edges.  Open edges are called \emph{parallel} if the corresponding rays
extend in the same direction.  Around a vertex, the counterclockwise order of
parallel edges is lost in the geometric representation, but encoded by the
permutation $\varsigma$.  

In this paper, we abuse the term \emph{graph} for the data $\Gamma=(\gH,
\involute, \rotate, \varrho)$, and we will also abuse the notation $\Gamma$ for
the image of $\varrho$.

\begin{remark}
	The structure can be modified into a graph represented in the Riemann sphere
	$\hat\C = \C \cup \{\infty\}$, as defined in~\cite{saddle1}, by adding an
	extra vertex at $\infty$ that closes all open edges.  This should help the
	readers to connect to the setup in~\cite{saddle1}.
\end{remark}

An orientation of the graph is a function $\sigma: \gH \to \{\pm 1\}$ such that
$\sigma(-h) = -\sigma(h)$ for all $h \in \gH \setminus \gR$.  A graph is said
to be \emph{orientable} if it has an orientation $\sigma$ such that
$\sigma\circ\rotate = -\sigma$.  In an orientable graph, every vertex has an
even degree.

\subsection{Discrete differential operators}

A (simple) cycle in the graph $\Gamma$ is a set of half-edges $c \subset \gH
\setminus \gR$ that can be ordered into a sequence $(h_1, \cdots, h_n)$ such
that $v(-h_i) = v(h_{i+1})$ for $1 \le i < n$ and $v(-h_n) = v(h_1)$, and
$v(h_i) \ne v(h_j)$ and $e(h_i) \ne e(h_j)$ whenever $i \ne j$.  The set of
cycles is denoted by $\gC$.  In particular, the orbits of $\rotate\involute$,
if contained in $\gH \setminus \gR$, are all cycles; we call these cycles
\emph{face cycles}, and use $\gF$ to denote the set of face cycles.  As the
graph is represented in the complex plane, we have necessarily
$|\gV|-|\gE|+|\gR|+|\gF|=1$.

A cut in $\Gamma$ is a set of half-edges $b \subset \gH$ such that, for some
fixed non-empty subset $\gV' \subseteq \gV$, we have $v(\eta) \in \gV'$ for all
$\eta \in b$ and $v(-h) \not\in \gV'$ for $h \in b \setminus \gR$.  The set of
cuts is denoted by $\gB$.  In particular, for any vertex $v$, the set
\[
	b(v) = \{ \eta \in v \colon v(-\eta) \ne v \quad \text{whenever} \quad \eta \notin \gR \}
\]
is a cut; we call these cuts \emph{vertex cuts}.

\begin{remark}[Notation]
	If $v(-h) = v(h)$, the edge $e(h)$ is a loop.  In this paper, graphs have no
	loops because they are represented in the complex plane $\C$.  So it makes
	sence to abuse the notation of vertex $v$ for the vertex cut $b(v)$.
\end{remark}

We use $\cA$ to denote the space of real-valued functions $f: \gH \setminus \gR
\to \R$ that are \emph{antisymmetric} in the sense that $f_{-h} = -f_h$ for all
$h \in \gH \setminus \gR$.  $\cA$ is a vector space of dimension $|\gE| -
|\gR|$.  Moreover, we use $\cR$ to denote the space of real-valued functions
$\gR \to \R$.

For $f \in \cA$, we define the discrete differential operator
\[
	\curl_c(f) = \sum_{h \in c} f_h, \qquad \curl(f)=(\curl_c(f))_{c \in \gC}.
\]
The image of $\curl$ is the \emph{cycle space} of $\Gamma$, and is denoted by
$\cC$.  It is a vector space of dimension $|\gF|$.  The projection $(x_c)_{c
\in \gC} \mapsto (x_c)_{c \in \gF}$ provides an isomorphism between $\cC$ and
$\R^{|\gF|}$.  In fact, the face cycles form a \emph{cycle basis};
see~\cite{saddle1}.

Let $f \in \cA \times \cR \simeq \R^{|\gE|}$ be a real-valued function on $\gH$
whose restriction on $\gH \setminus \gR$ is antisymmetric.  We define the
operator
\[
	\div_b(f) = \sum_{h \in b} f_h, \qquad \div(f)=( \div_b(f))_{b \in \gB}.
\]
The image of $\div$ is the \emph{cut space} of $\Gamma$, and is denoted by
$\cB$. It is a vector space of dimension $|\gV|$.  The projection $(x_b)_{b \in
\gB} \mapsto (x_v)_{v \in \gV}$ provides an isomorphism between $\cB$ and
$\R^{|\gV|}$.  In fact, the vertex cuts form a \emph{cut basis};
see~\cite{saddle1}.

For each half-edge $h \in \gH \setminus \gR$, let $\cent\ell_h$ be the length
of the segment $\varrho(e(h))$.  For $b \in \gB$, define
\[
	\cent\ell_b= \min_{h \in b \setminus \gR} \cent\ell_h\quad\text{and}\quad m(b)=\{ h \in b \setminus \gR \mid \cent\ell_h = \cent\ell_b\}.
\]
For $\phi \in \cA$, we define the operator
\[
	\mdiv_b(\phi)=\sum_{h\in m(b)}\phi_h, \qquad
	\mdiv(\phi) = (\mdiv_b(\phi))_{b\in\gB}.
\]
The same argument as in~\cite{saddle1} proves that the image of $\mdiv$,
denoted by $\cB_m$, has the same dimension as $\cB$.  In particular, there is a
cut basis $\gB_m^* \subset \gB$ such that the projection $(x_b)_{b \in \gB}
\mapsto (x_b)_{b \in \gB_m^*}$ provides an isomorphism between $\cB_m$ and
$\R^{|\gV|}$.

\subsection{Horizontal balance and rigidity}

To each $h \in \gH \setminus \gR$, we assign the unit tangent vector
$\cent{u}_h=e^{\ii\cent\theta_h}$ of the segment $\varrho(e(h))$ at
$\varrho(v(h))$. We denote by $\cent\ell_h$ the length of the segment
$\varrho(e(h))$ and set $\cent{x}_h=\cent\ell_h\cent{u}_h$.

For a ray $r \in \gR$, $\cent x_r$ and $\cent\ell_r$ are not defined.  It is
only assigned a unit vectors $\cent u_r = e^{\ii\cent\theta_r}$ in the
direction of the ray $\varrho(e(r))$.

\begin{remark}[Notation]
	We distinguish the notations $\vartheta=(\theta_r)_{r \in \gR}$ and
	$\theta=(\theta_\eta)_{\eta \in \gH}$.  They are both frequently used in this
	paper.
\end{remark}

For $\chi = (x, \vartheta) \in \cA^2 \times \cR \simeq \R^{|\gH|}$ in a
neighborhood of $\cent \chi = (\cent x, \cent\vartheta)$, we define
\[
	u_\eta(\chi)= \begin{dcases}
		\frac{x_\eta}{\|x_\eta\|}, & \eta \in \gH \setminus \gR,\\
		\exp(\ii\theta_\eta), & \eta \in \gR,
	\end{dcases}
	\quad \text{and}\quad
	u(\chi) = (u_\eta(\chi))_{\eta \in \gH}.
\]

The horizontal periods are given by the function
\begin{align*}
	\Phor\colon \cA^2 \times \cR & \to \cC^2 \\
	\chi = (x, \vartheta) & \mapsto \curl(x).
\end{align*}
As the graph $\Gamma$ is represented in the complex plane $\C$, $\cent\chi$
solves the \emph{horizontal period problem}.
\[
	\Phor_c(\cent\chi) = 0,\quad c \in \gF.
\]

The horizontal forces are given by the function
\begin{align*}
	\Fhor\colon \cA^2 \times \cR & \to \cB^2 \\
	\chi = (x,\vartheta) & \mapsto \div(u(\chi)).
\end{align*}

\begin{definition}
	The graph $\Gamma$ is balanced if $\Fhor(\cent\chi) = 0$, and is rigid if
	\[
		(D\Fhor(\cent\chi),\Phor) \colon \cA^2 \times \cR \to \cB^2 \times \cC^2
		\simeq \cA^2 \times \R^2
	\]
	is surjective.
\end{definition}

If the graph is balanced and rigid, then in the neighborhood of $\cent\chi$,
the set of $\chi$ that solves $(\Fhor, \Phor)=0$ form a manifold $\cM$ of
dimension~$(|\gR|-2)$.

In an orientable and balanced graph, we say that a vertex $v$ is
\emph{degenerate} if the unit vectors $\cent{u}_h$, $h \in v$, are collinear;
we say that $v$ is \emph{special} if $\deg(v) \ge 6$ and $\deg(v)-2$ of the
unit vectors $\cent{u}_h$ are collinear.  We say that $v$ is ordinary if it is
neither degenerate nor special.

We want to place a saddle tower $\sS_v$ at each vertex $v \in \gV$ with their
wings along the edges in $\gE$.  This is possible only if the graph is
orientable, balanced, and all vertices are ordinary.  Then the following
proposition, whose proof is delayed to the appendix, asserts that the
horizontal rigidity is guaranteed in the absence of parallel edges.

\begin{proposition}\label{prop:rigid}
	If the graph $\Gamma$ is orientable, balanced, has no parallel edges, and all
	vertices are ordinary, then $\Gamma$ is rigid.
\end{proposition}

\begin{remark}
	A graph with parallel edges might not be rigid; see Example~\ref{ex:benzene}.
\end{remark}

The phase of a saddle tower, informally speaking, is the height of its
horizontal reflection plane; we recommend the readers to~\cite{saddle1} for the
formal definition.  The phase differences between saddle towers are prescribed
through an antisymmetric phase function $\cA \ni \cent\phi \colon \gH \setminus
\gR \to \R/2\pi\Z$.  We say that $\cent\phi$ is \emph{trivial} if $\cent\phi =
0$ or $\pi$ on every half-edge.  Trivial phase functions give rise to SPMSs
with horizontal symmetry planes, as claimed in the following theorem.

\begin{theorem}[SPMSs with horizontal symmetry plane]\label{thm:symmetry}
 	Given a graph $\Gamma$ and a trivial phase function $\cent\phi$.  Assume that
 	$\Gamma$ is orientable, balanced, rigid, and all vertices are ordinary.  Then
 	for sufficiently small $\varepsilon$, there is a continuous family
 	$\sM_\varepsilon$ of immersed singly periodic minimal surfaces of genus
 	$|\gF|$ in $\R^3$ with $|\gR|$ Scherk ends, vertical period $(0,0,2\pi)$, and
 	a horizontal symmetry plane, such that
 	\begin{enumerate}
		\item $\varepsilon^2\sM_\varepsilon$ (scaling of $\sM_\varepsilon$ by
			$\varepsilon^2$) converges to $\Gamma \times \R$ as $\varepsilon\to 0$.

		\item For each vertex $v \in \gV$, there exists a horizontal vector
			$X_v(\varepsilon)$ such that $\sM_\varepsilon-X_v(\varepsilon)$
			converges on compact subsets of $\R^3$ to a saddle tower $\sS_v$ as
			$\varepsilon \to 0$. Moreover, $\varepsilon^2X_v(\varepsilon) \to
			\varrho(v)$ as $\varepsilon\to 0$.

		\item For each half-edge $h \in \gH \times \gR$, the phase difference of
			$\sS_{v(-h)}$ over $\sS_{v(h)}$ is equal to $\cent\phi_h$.
	\end{enumerate}
\end{theorem}

In fact, this family also depend continuously on $\chi \in \cM$ in a
neighborhood of $\cent\chi$.

\subsection{Vertical balance and rigidity}

In the following, we will prescribe non-trivial phase functions $\cent\phi$.
We define the vertical periods as the function
\begin{align*}
	\Pver\colon \cA &\to \cC\\
	\phi &\mapsto \curl(\phi).
\end{align*}
We require that the prescribed phase function~$\cent\phi$ solve the
\emph{vertical period problem}
\[
	\Pver_c(\cent\phi) = 0, \quad c \in \gF.
\]

We now explain the vertical balancing condition.  For each vertex $v \in \gV$,
consider a punctured Riemann sphere $\C_v$ on which the Weierstrass
parameterization of $\sS_v$ is defined.  Recall from~\cite{saddle1} that the
punctures must lie on a circle fixed by the anti-holomorphic involution $\rho$
corresponding to the reflection symmetries of $\sS_v$ in horizontal planes.
Then for each $\eta \in v$, fix a local coordinate $w_\eta$ in a neighborhood
of the puncture $p_\eta \in \C_v$ vanishing at $p_\eta$, and assume that
$w_\eta$ is \emph{adapted} in the sense that $w_\eta \circ \rho =
\overline{w_\eta}$.  Recall from~\cite{saddle1} the quantities $\Upsilon_\eta$
are $\mu_\eta$ that describe the shape of the wings; they are defined in terms
of the Weierstrass parameterization and the local coordinates $w_\eta$.

In the following, we write $\cD := \ker(D\Fhor(\cent\chi), \Phor)$ for the
space of infinitesimal deformations of the graph that preserve the balance.  If
the graph is rigid, then $\cD$ is a subspace of dimension $|\gR|-2$.  Note that
$\cD$ include the rotations and scalings of the graph.

Let $\dot\zeta$ be the unique solution in $\cD^\perp$ (with respect to the
standard inner product of $\cA^2 \times \cR \simeq \R^{|\gH|+|\gR|}$) to the
linear system
\begin{equation}\label{eq:system}
	\begin{cases}
		\Phor_c(\dot\zeta) = - \Phor_c(\mu^a), & c \in \gF,\\
		D\Fhor_v(\cent\chi) \cdot \dot\zeta = 0, & v \in \gV,
	\end{cases}
\end{equation}
where $\mu^a_h = \mu_h - \mu_{-h}$.  Then a general solution
to~\eqref{eq:system} is of the form $\xi + \dot\zeta$, where $\xi \in \cD$.
Fix a prescribed $\cent\xi \in \cD$, we write $\cent\xi + \dot\zeta =:
\cent{\dot\chi} = (\cent{\dot x}, \cent{\dot\vartheta})$, and define
\[
	K_h = \Upsilon_h \Upsilon_{-h} e^{-\re(\cent{\dot x}_h\overline{\cent{u}_h})}
\]
for $h \in \gH \setminus \gR$.  As in~\cite{saddle1}, $K_h$ is real, positive,
independent of horizontal translations of the saddle towers and the adapted
local coordinates $w_h$.  It depends on $\cent\xi$, but the dependence is
omitted for simplicity.  When their concrete values matter but are not
specified in the context, it is understood that $\cent\xi = 0$.

We define the vertical forces as the function
\begin{align*}
	\Fver \colon \cA & \to \cB_m \\
	\phi &\mapsto \mdiv\big( (K_h \sin\phi_h)_{h \in \gH \setminus \gR} \big).
\end{align*}

\begin{definition}
	The phase function $\cent\phi$ is \emph{balanced} if $\Fver(\cent\phi) = 0$,
	and is \emph{rigid} if $(D\Fver(\cent\phi),\Pver)$ is an isomorphism.
\end{definition}

\begin{remark}\label{rmk:continuous}
	In general, the vertical forces do not depend continuously on the parameter
	$\cent\chi$ that describes the graph, but they depend continuously on the
	infinitesimal deformation $\cent\xi$.
\end{remark}

\begin{remark}\label{rmk:normalize2}
	The vertical balance is invariant under the transformation
	\[
		(\cent{\dot x}, \cent{\dot\vartheta}) \mapsto
		(\cent{\dot x} + \lambda \cent x, \cent{\dot\vartheta} + \arg\lambda)
	\]
	for $\lambda \in \C$.
\end{remark}

The following proposition adapted from~\cite{saddle1} shows that vertical
balance is not rare.  A similar proof as in~\cite{saddle1} applies here with
little modification.

\begin{proposition} \label{proposition:zero-is-rigid}
	Let $\cent\phi$ be a phase function such that $\cos(\cent\phi_h)$, $h \in \gH
	\setminus \gR$, are all positive or all negative.  Then $\cent\phi$ is rigid.
	In particular, the zero phase function is always balanced and rigid.
\end{proposition}

We call the pair $(\Gamma, \cent\phi)$ a \emph{configuration} and say that the
configuration is balanced (resp.\ rigid) if both $\Gamma$ and $\cent\phi$ are
balanced (resp.\ rigid).  Our main result for SPMSs is the following.

\begin{theorem}[SPMSs]\label{thm:main}
 	Given a configuration $(\Gamma, \cent\phi)$ and a prescribed deformation
 	$\cent\xi \in \cD$.  Assume that $\Gamma$ is orientable, the configuration is
 	balanced and rigid, and all vertices are ordinary.  Then for
 	$(\varepsilon,\xi) \in \R_+ \times \cD$ in a neighborhood of $(0,\cent\xi)$,
 	there is a continuous family $\sM_{\varepsilon,\xi}$ of immersed singly
 	periodic minimal surfaces of genus $|\gF|$ in $\R^3$ with $|\gR|$ Scherk ends
 	and vertical period $(0,0,2\pi)$ such that
 	\begin{enumerate}
		\item $\varepsilon^2\sM_{\varepsilon,\xi}$ (scaling of
			$\sM_{\varepsilon,\xi}$ by $\varepsilon^2$) converges to $\Gamma\times\R$
			as $\varepsilon\to 0$.  In particular, the $|\gR|$ Scherk ends of
			$\varepsilon^2\sM_{\varepsilon,\xi}$ converge to the rays of $\Gamma$.

		\item For each vertex $v \in \gV$, there exists a horizontal vector
			$X_v(\varepsilon)$ such that $\sM_{\varepsilon,\xi}-X_v(\varepsilon)$
			converges on compact subsets of $\R^3$ to a saddle tower $\sS_v$ as
			$\varepsilon \to 0$. Moreover, $\varepsilon^2X_v(\varepsilon) \to
			\varrho(v)$ as $\varepsilon\to 0$.

		\item For each half-edge $h \in \gH \setminus \gR$, the phase difference of
			$\sS_{v(-h)}$ over $\sS_{v(h)}$ is equal to $\phi_h$.
	\end{enumerate}
\end{theorem}

\begin{remark}
	In view of Remarks~\ref{rmk:normalize2}, we may, up to scalings and
	horizontal rotations, fix $\dot x_h = 0$ for a particular half-edge $h \in \gH
	\setminus \gR$.  Then we construct families with $|\gR|-3$ parameters.
	Recall that the Karcher-Scherk saddle towers with $2k$ Scherk ends form a
	family with $2k-3$ parameters.
\end{remark}

\subsection{Embeddedness}

We identify several scenarios where the surfaces in Theorem~\ref{thm:main} are
embedded.

In the first case, we assume that the graph has no parallel rays, hence the
Scherk ends of the surfaces do not intersect for sufficiently small
$\varepsilon$.  This is the case considered in~\cite{traizet2001}.

For a more formal statement, let us label the rays by integers $r=1, \cdots,
|\gR|$, in the counterclockwise order.  Recall that, around the same vertex, the
counterclockwise order of parallel edges is given by the permutation $\varsigma$.
Up to a rotation, we may assume that $0 = \cent\theta_1 \le \cdots \le
\cent\theta_{|\gR|} < 2\pi$.  Two rays $r$ and $r'$ are then parallel if
$\cent\theta_r = \cent\theta_{r'}$.

\begin{theorem}\label{thm:embed1}
	The minimal surfaces $M_{\varepsilon,\xi}$ in Theorem~\ref{thm:main} is
	embedded for $(\varepsilon,\xi)$ sufficiently close to $(0,\cent\xi)$ if
	$\cent\theta_r < \cent\theta_{r+1}$ for all $1 \le r < |\gR|$.
\end{theorem}

If the graph has parallel rays, the corresponding Scherk ends become parallel
in the limit $\varepsilon \to 0$.  We call them \emph{parallel Scherk ends}.  

\begin{remark}
	This terminology is certainly an abuse.  As we have stressed several times,
	the Scherk ends might bend from the direction of the corresponding rays,
	hence might not be parallel for $\varepsilon > 0$!
\end{remark}

We want to \emph{resolve} these parallel Scherk ends.  That is, as
$\varepsilon$ increases, we want the Scherk ends to bend away from each other.
Recall that the lowest order deformation of the graph is prescribed by
$\cent{\dot\chi}= (\cent{\dot x}, \cent{\dot\vartheta})$.

\begin{theorem}\label{thm:embed2}
	The minimal surfaces $M_{\varepsilon,\xi}$ in Theorem~\ref{thm:main} is
	embedded for $(\varepsilon,\xi)$ sufficiently close to $(0,\cent\xi)$ if
	$\cent{\dot\theta}_r < \cent{\dot\theta}_{r+1}$ whenever $\cent\theta_r =
	\cent\theta_{r+1}$, $1 \le r < |\gR|$.
\end{theorem}

If parallel Scherk ends are not resolved by the lowest order deformation of the
graph, it might still be resolved by higher order terms in the Taylor expansion
of $\chi$ in $\varepsilon$.  Eventually, we may use the entire Taylor series of
$\chi$.  More formally, fix an analytic function $\xi(\varepsilon) \in \cD$
such that $\xi(0) = \cent\xi$.  Since the graph is balanced and rigid, the
non-linear system
\begin{equation}\label{eq:systemanal}
	\begin{cases}
		\Phor(\widetilde\zeta(\varepsilon)) = - \varepsilon^2 \Phor(\mu^a),\\
		\Fhor(\widetilde\chi(\varepsilon)) = 0,
	\end{cases}
\end{equation}
where $\widetilde\chi(\varepsilon) = \cent\chi + \varepsilon^2 \xi(\varepsilon) +
\widetilde\zeta(\varepsilon) = (\widetilde x(\varepsilon),
\widetilde\vartheta(\varepsilon))$ is the limit of the Taylor series of $\chi$,
has a unique solution $\widetilde\zeta(\varepsilon) \in \cD^\perp$ for
$\varepsilon$ sufficiently small.

\begin{theorem}\label{thm:embed3}
	The $1$-parameter family of minimal surfaces
	$M_{\varepsilon,\xi(\varepsilon)}$ as constructed in Theorem~\ref{thm:main}
	is embedded for sufficiently small $\varepsilon$ if
	$\widetilde\theta_r(\varepsilon) < \widetilde\theta_{r+1}(\varepsilon)$ for
	sufficiently small $\varepsilon$ whenever $\cent\theta_r =
	\cent\theta_{r+1}$, $1 \le r < |\gR|$.
\end{theorem}

Interestingly, if the graph appears as a simple line arrangement, as considered
in~\cite{traizet2001}, then even the Taylor series $\widetilde\chi$ is not
enough to resolve parallel Scherk ends; see Lemma~\ref{lem:lines}.  In this
case, we must consider deformations that are not analytic but flat in
$\varepsilon$.  As $\varepsilon$ increases, the horizontal forces will slightly
deviate from unit vectors.  For the force along the half-edge $h$, the
deviation is in the order of $\exp(-\ell_h/\varepsilon^2)$.  The surface must
deform to balance the horizontal forces, and this deformation might resolve
parallel Scherk ends.  See, for instance, Example~\ref{ex:3line}.

More formally, define
\[
	m(\gH)=\{ h \in \gH \setminus \gR \mid \cent\ell_h \le \cent\ell_{h'}\, \forall h' \in \gH \setminus \gR \}.
\]
Let $\widehat\zeta$ be the unique solution in $\cD^\perp$ to the linear system
\begin{equation}\label{eq:systemflat}
	\begin{dcases}
		\Phor_c(\widehat\zeta) = - \sum_{h \in c \cap m(\gH)} x_h K_h \cos \cent\phi_h, & c \in \gF,\\
		D\Fhor_v(\cent\chi) \cdot \widehat\zeta = -\sum_{h \in v \cap m(\gH)} u_h K_h \cos \cent\phi_h, & v \in \gV,
	\end{dcases}
\end{equation}
and write $\widehat\zeta = (\widehat x, \widehat\vartheta)$.

\begin{theorem}\label{thm:embed4}
	If the graph $\Gamma$ appears as a simple line arrangement, then the
	$1$-parameter family of minimal surfaces $M_{\varepsilon,\xi(\varepsilon)}$
	as constructed in Theorem~\ref{thm:main} is embedded for sufficiently small
	$\varepsilon$ if, whenever $\cent\theta_r = \cent\theta_{r+1}$, we have
	$\widetilde\theta_r(\varepsilon) = \widetilde\theta_{r+1}(\varepsilon)$ for
	sufficiently small $\varepsilon$ and $\widehat\theta_r <
	\widehat\theta_{r+1}$.
\end{theorem}

There are certainly situations where even Theorem~\ref{thm:embed4} can not
determine embeddedness; see Example~\ref{ex:verysubtle} below.  If this is the
case, parallel Scherk ends might still be resolved by looking at even higher
orders, but we do not plan to continue. 

\section{Examples}\label{sec:examples}

\subsection{Trees}\label{ssec:tree}

We say that $\Gamma$ is a \emph{tree} if it has no cycle of length $>2$.  Note
that we allow cycles of length $2$ in the trees; edges in such a cycle must be
parallel closed edges.

\begin{lemma}\label{lem:tree}
	If $\Gamma$ is a tree, then the only balanced phase functions on $\Gamma$ are
	the trivial ones.
\end{lemma}

\begin{proof}
	Recall that parallelism is an equivalence relation on the set of closed
	edges.  For edges in a parallelism class, $\cent\phi$ must take the same
	value on their half-edges from the same vertex, so that the vertical period
	problem is solved. If $\Gamma$ is a tree, then these half-edges form a cut.
	So a balanced $\cent\phi$ must be $0$ or $\pi$ on every half-edge.
\end{proof}

Then by Theorem~\ref{thm:symmetry}, if a tree is orientable, balanced, rigid,
and has only ordinary vertices, it always gives rise to symmetric SPMSs.  We
see here that the symmetry is induced by the structure.  If the tree has no
cycle of length $2$, it gives rise to a SPMS with genus zero and Scherk ends.
In this case, it was proved in~\cite{perez2007} that the symmetry is imposed by
the structure.

\medskip

\addtocounter{example}{-1}

\begin{example}[revisit]
	The graph in Figure~\ref{fig:tree1} is a tree, so the phase difference
	$\cent\phi$ between the two saddle towers is either $0$ (in phase) or $\pi$
	(in opposite phase).

	The graph is given by a simple line arrangement.  By Lemma~\ref{lem:lines}
	below, under the analytic deformation $\widetilde\chi$ that
	solves~\eqref{eq:systemanal}, the directions of two rays on the same line
	must remain opposite.  As a consequence, if the two Scherk ends pointing
	upwards bend away from each other under $\widetilde\chi$, the Scherk ends
	pointing downwards must bend towards each other, creating unwanted
	intersection.  The only way to avoid this is to let the vertical rays remain
	parallel under $\widetilde\chi$.

	Without loss of generality, we fix the analytic function $\xi(\varepsilon)
	\equiv 0$.  This guarantees that the vertical rays remain vertical under
	$\widetilde\chi$, and the horizontal rays remain horizontal.  It remains to
	study non-analytic deformation by solving~\eqref{eq:systemflat}.  In this
	simple example, there is no equation for $\Phor_c$, and $m(\gH)$ consists of
	only half-edges in the middle edge.  One easily solves a system of two
	equations and finds that the solution $\widehat\zeta \in \cD^\perp$,
	depending on the phase difference $\cent\phi$, is as depicted by dashed
	graphs in Figure~\ref{fig:tree1}.

	The solution can also be understood ``physically''.  We will see that, as
	$\varepsilon$ increases, the force along the middle edge varies by
	$K_h\cos(\cent\phi)\exp(-\cent\ell/\varepsilon^2)$, where $K_h > 0$ and
	$\cent\ell$ is the length of the edge.  If $\cent\phi = 0$, the force
	increases, while the forces along the rays remain unit vectors.  So the
	vertical Scherk ends must bend away from each other to recover horizontal
	balance, and the surfaces are embedded by Theorem~\ref{thm:embed4}.
	Otherwise, if $\cent\phi = \pi$, then the force along the middle edge
	decreases, and the vertical Scherk ends must bend towards each other to
	recover horizontal balance, and the surfaces are not embedded.  In no case do
	the Scherk ends remain vertical, contradicting~\cite{morabito2020}.  \qed
\end{example}

\begin{figure}[!htb]
	\includegraphics[width=.4\textwidth]{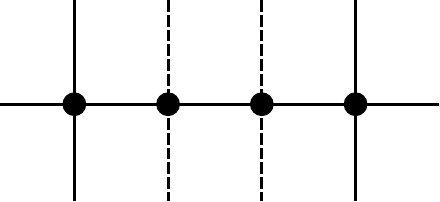}
	\caption{
		A configuration whose embeddedness can not be determined even by
		Theorem~\ref{thm:embed4}. \label{fig:tree2}
	}
\end{figure}

\begin{example}[A limitation]\label{ex:verysubtle}
	Figure~\ref{fig:tree2} shows an example that only adds two more vertical
	lines (dashed) to Example~\ref{ex:3line}.  Assume that $\cent\phi = 0$ on all
	half-edges, and that $\xi(\varepsilon) \equiv 0$.  This is a situation that
	even Theorem~\ref{thm:embed4} can not determine embeddedness.  To see this,
	note that the horizontal forces increase, to the lowest order, by the same
	amount on all closed edges.  As a consequence, the Scherk ends corresponding
	to solid vertical rays can not remain vertical as claimed
	in~\cite{morabito2020}, but must bend outwards to balance the horizontal
	forces.  However, the Scherk ends corresponding to dashed vertical rays must
	stay vertical.  We can not tell how the middle Scherk ends bend, or if they
	bend at all, without looking at higher-order terms.  But we do not plan to do
	this. \qed
\end{example}

\begin{figure}[!htb]
	\includegraphics[width=.4\textwidth]{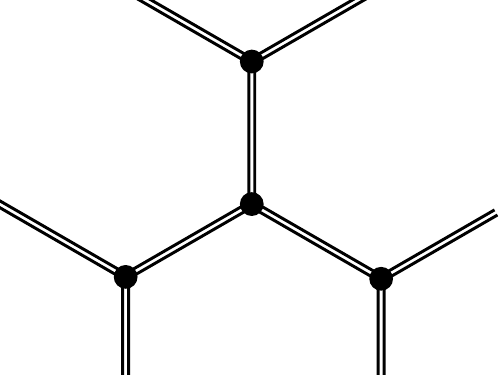}
	\caption{A tree with parallel edges.  \label{fig:tree3}}
\end{figure}

\begin{example}[Another limitation]\label{ex:paralleltree}
	Figure~\ref{fig:tree3} shows an example with parallel edges.  It does not
	appear as a line arrangement, but suffers the same problem regarding
	embeddedness: The parallel Scherk ends can not be resolved by the analytic
	deformation $\widetilde\chi$.  If $\xi(\varepsilon) \equiv 0$ and the phase
	function $\cent\phi$ is $0$ on all half-edges, then the forces along closed
	edges will all increase with $\varepsilon$.  To recover horizontal balance,
	parallel rays from different vertices must bend away from each other.
	
	But another problem arises: We can not tell how parallel rays from the same
	vertex bend.  To study this behavior, one needs to look closer into the
	detailed structure of the 6-wing Karcher-Scherk saddle towers.  We do not
	plan to do this in the current manuscript.  \qed
\end{example}

\subsection{Previously known examples}

The phase function $\cent\phi$ can be recovered from a phase function
$\cent\varphi: \gV \to \R/2\pi\Z$, unique up to the addition of a constant,
such that $\cent\phi_h = \varphi_{v(-h)} - \varphi_{v(h)} \pmod{2\pi}$.  In
this paper, the phase functions are marked in the figures by labelling
$\cent\varphi_v$ on the vertices.

In this part, we show some classical examples that may arise from our
construction.

\begin{figure}[!htb]
	\includegraphics[width=.3\textwidth]{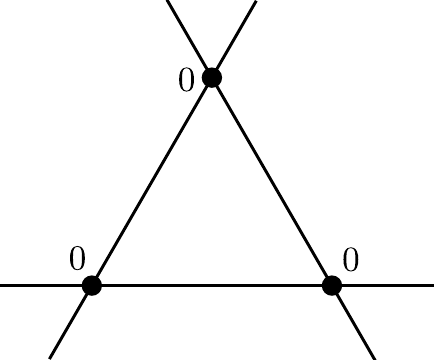}
	\hspace{0.5cm}
	\includegraphics[width=.3\textwidth]{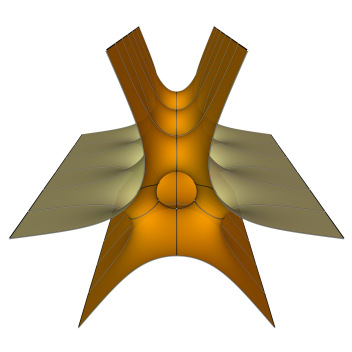}

	\includegraphics[width=.3\textwidth]{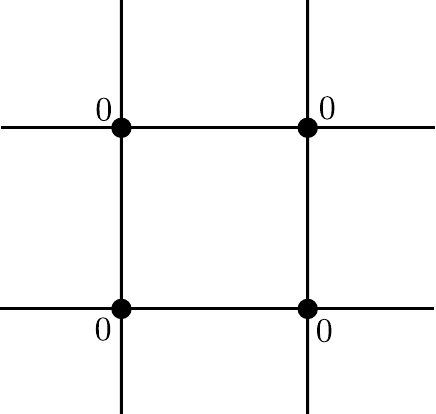}
	\hspace{0.5cm}
	\includegraphics[width=.3\textwidth]{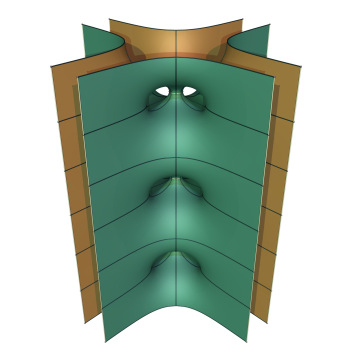}

	\caption{Karcher--Scherk saddle towers with handles. \label{fig:handle}}
\end{figure}

\begin{example}[Toroidal saddle towers]\label{ex:handle}
	Figure~\ref{fig:handle} shows symmetric Karcher--Scherk saddle towers with a
	vertical tunnel in the middle.  They were first constructed by
	Karcher~\cite{karcher1988} with explicit Weierstrass data defined on tori,
	hence the name.  In the framework of our construction, the one with 6 ends
	arises from three lines forming an equiangular triangle, the one with 8 ends
	arises from four lines forming a square.

	Note that the one with 8 ends has parallel ends, but they do not correspond
	to the parallel rays in the graph.  In fact, we have $\xi(\varepsilon)\equiv 0$ and
	the Scherk ends bend from vertical and horizontal directions in the limit
	$\varepsilon \to 0$ all the way to the diagonal directions as $\varepsilon$
	increases, forming new parallel pairs.  If $\varepsilon$ increase
	further, the Scherk ends will intersect.  \qed
\end{example}

\begin{figure}[!htb]
	\includegraphics[width=.3\textwidth]{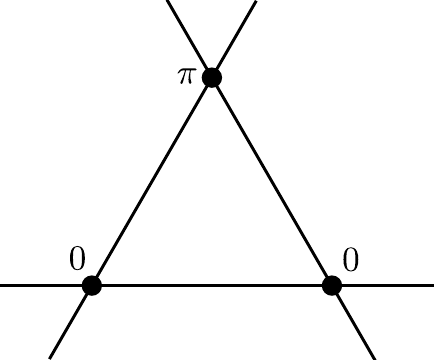}
	\hspace{0.5cm}
	\includegraphics[width=.3\textwidth]{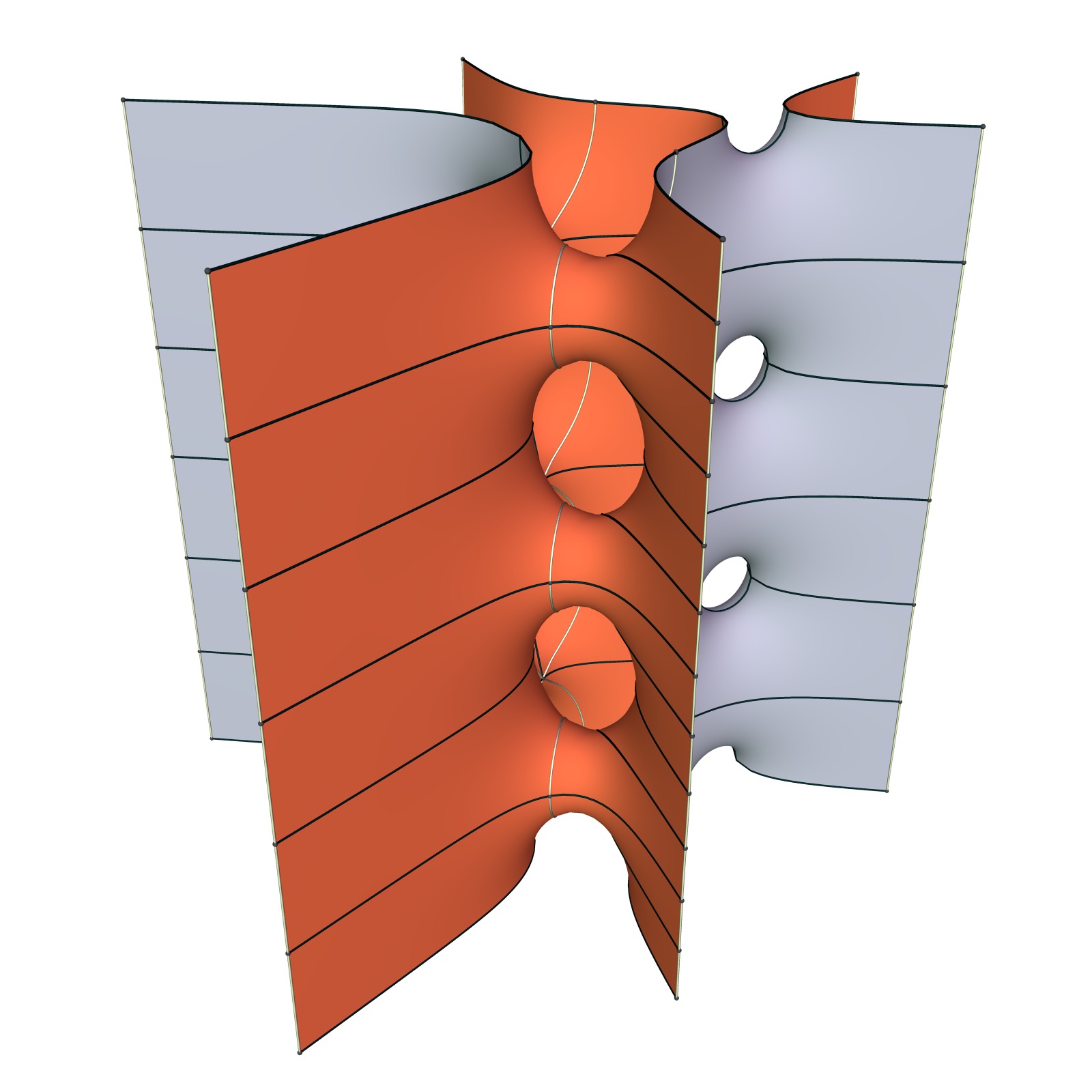}

	\caption{The Costa--Scherk surfaces.\label{fig:CS}}
\end{figure}

\begin{figure}[!htb]
	\includegraphics[width=.25\textwidth]{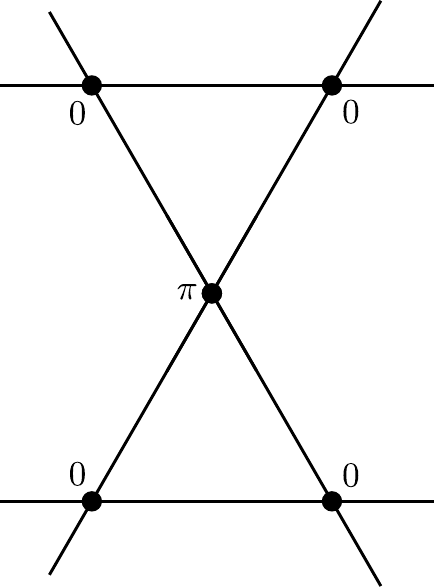}
	\hspace{0.5cm}
	\includegraphics[width=.35\textwidth]{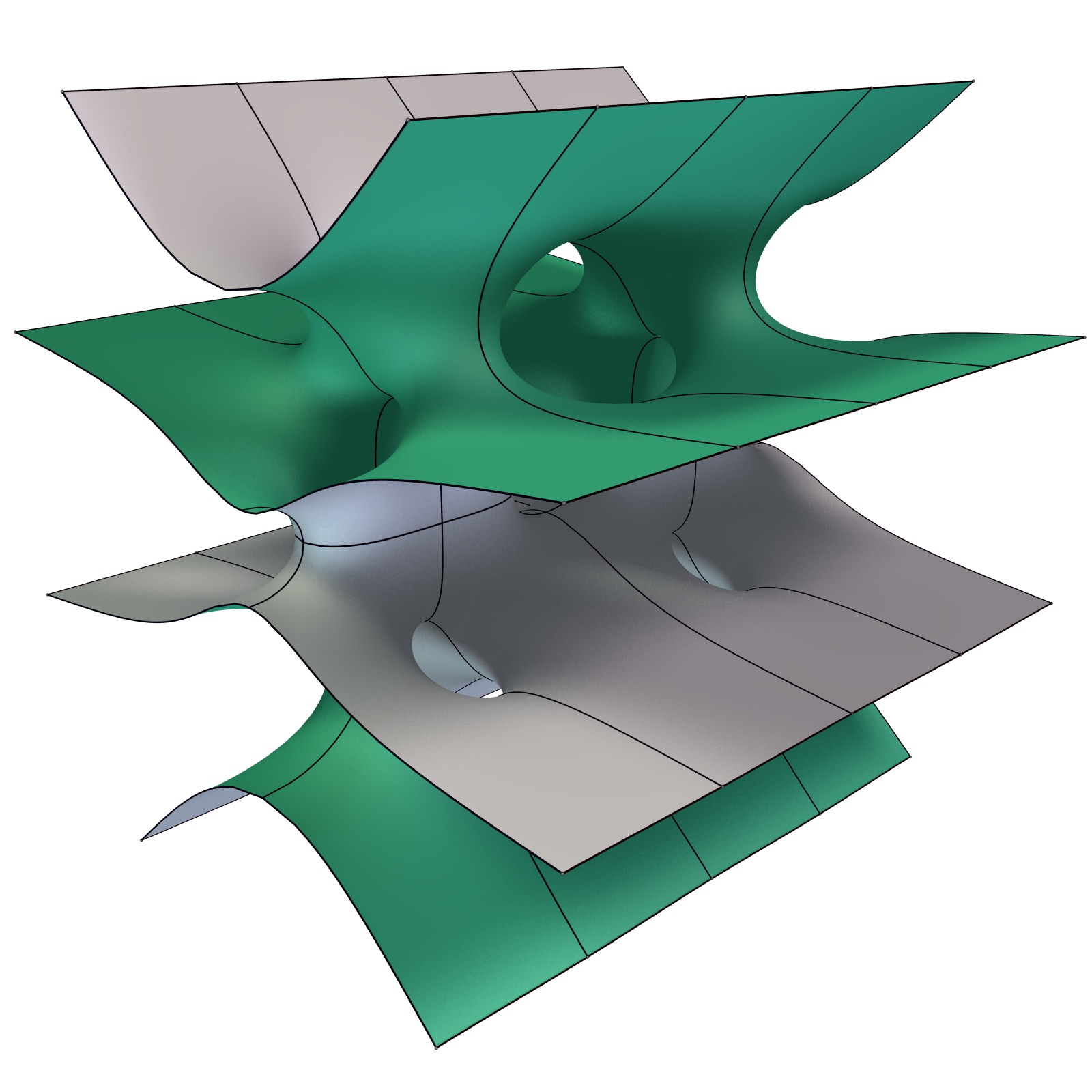}

	\caption{The da Silva--Batista surface.\label{fig:dSB}}
\end{figure}

\begin{example}[Costa--Scherk surfaces]
	In Figure~\ref{fig:CS} is a family that arises again from three lines forming
	an equiangular triangle.  But this time we use a different trivial phase
	function.  The SPMSs looks like a tower of Costa surfaces, hence the name.
	\qed
\end{example}

\begin{example}[Da Silva--Batista surfaces]
	In Figure~\ref{fig:dSB} is a family that arises from an arrangement of four
	lines.  As we have 8 ends, the family is described by 5 parameters up to
	scalings and Euclidean motions.  In particular, it includes the saddle tower
	limit of the 2-parameter family constructed in~\cite{dasilva2010}, for which
	vertical reflection planes were assumed.  \qed
\end{example}

\subsection{Singly periodic gyroids}\label{ssec:gyroid}

\begin{figure}[!htb]
	\includegraphics[width=.3\textwidth]{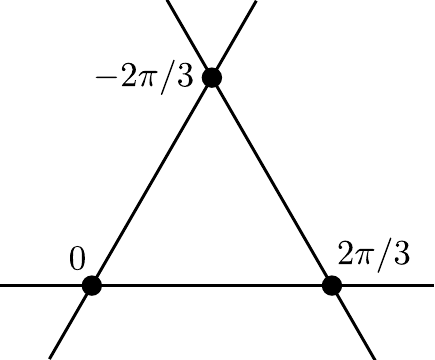}
	\hspace{0.5cm}
	\includegraphics[width=.3\textwidth]{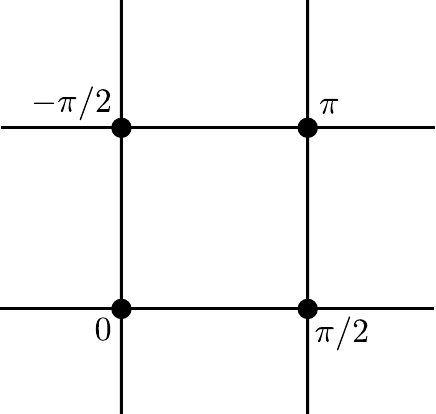}
	\caption{Singly periodic analogues of rGL and tG surfaces.\label{fig:singlyG}}
\end{figure}

\begin{example}[Singly periodic rGL]
	In view of Lemma~\ref{lem:tree}, a SPMS with Scherk ends and no horizontal
	symmetry plane must arise from a graph with a cycle of length at least three.
	The configuration on the left of Figure~\ref{fig:singlyG}, where three lines
	form an equiangular triangle, is therefore the smallest non-symmetric,
	balanced, and rigid example.  It can be seen as the singly periodic analogue
	of the rGL family of triply periodic minimal surfaces; see~\cite{chen2019}
	and~\cite{saddle1}.  Note that, if the triangle was not equiangular, then the
	only possible balanced phase functions are the trivial ones, and we obtain a
	deformation of the toroial saddle tower with 6 ends (see
	Example~\ref{ex:handle}).  One sees here that the vertical balance does not
	depend continuously on the graph; see Remark~\ref{rmk:continuous}. \qed
\end{example}

\begin{example}[Singly periodic tG, inconclusive]
	In the same spirit, the configuration on the right of
	Figure~\ref{fig:singlyG}, where four lines form a square, can be seen as the
	singly periodic analogue of the tG family of triply periodic minimal
	surfaces; see~\cite{chen2019} and~\cite{saddle1}.  This configuration is
	balanced, but not vertically rigid: It seems that one may vertically slide
	two non-adjacent saddle towers with respect to the others without any
	horizontal deformation of the graph.  Hence our construction is not
	conclusive on this configuration.  Even if this configuration does give rise
	to SPMSs with Scherk ends, it would still be challenging to determine their
	embeddedness, as Theorem~\ref{thm:embed4} does not apply here. \qed
\end{example}

\subsection{Miscellaneous examples}

\begin{figure}[!htb]
	\includegraphics[width=.4\textwidth]{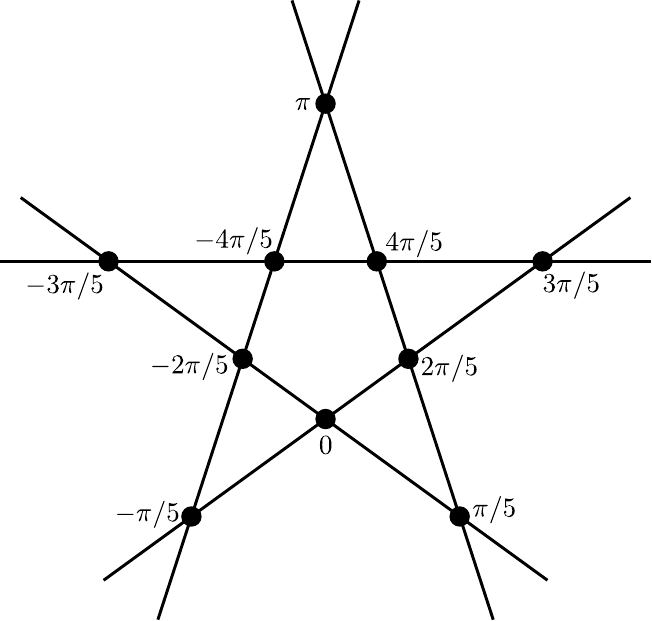}
	\hspace{0.5cm}
	\includegraphics[width=.5\textwidth]{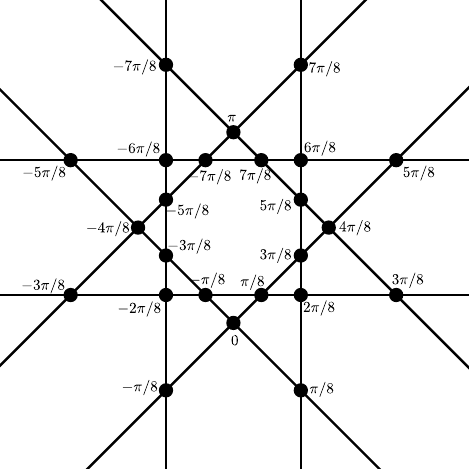}
	\caption{Non-symmetric examples from polygrams.\label{fig:polygram}}
\end{figure}

\begin{example}[Polygrams]
	Non-symmetric examples can be produced from the graph that appears as a
	simple arrangement of $k$ lines, $k > 4$, that form a regular polygram.  More
	specifically, such a graph contains a clockwise cycle whose edges are all of
	the shortest length.  We choose $\cent\phi$ that takes the same value on all
	half-edges in this cycle, and its value on other half-edges can be determined
	(not necessarily unique!) by solving the period and balance problems.
	Figure~\ref{fig:polygram} shows two examples with $k=5$ and $k=8$.  These
	configurations all appear as line arrangements.  If $k$ is odd, the
	embeddedness follows from Theorem~\ref{thm:embed1}.  If $k$ is even, and
	$\xi(\varepsilon) \equiv 0$, the embeddedness can be determined by
	Theorem~\ref{thm:embed4}.  \qed
\end{example}

\begin{figure}[!htb]
	\includegraphics[width=.4\textwidth]{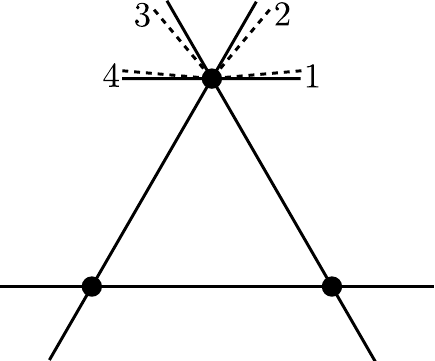}
	\caption{
		An infinitesimal deformation $\cent\xi$ that resolve parallel Scherk ends.
		\label{fig:misc1}
	}
\end{figure}

\begin{example}
	Figure~\ref{fig:misc1} shows a graph with two vertices of degree 4 and one
	vertex of degree 6.  Let $c$ be the unique counterclockwise cycle in the graph.
	Using the explicit values of $\mu_h$ computed in~\cite{saddle1}, one verifies
	that $\Phor_c(\mu^a) = 0$.  So we have, very conveniently, $\dot\zeta = 0$
	and $\cent{\dot \chi} = \cent\xi \in \cD$.

	Label four rays as shown in the Figure.  One verifies that $\xi = (\dot x,
	\dot \vartheta)$ with $\dot x = 0$,
	\[
		2\dot\theta_1 = -\dot\theta_2
		= \dot\theta_3 = -2 \dot\theta_4 > 0,
	\]
	and $\dot\theta_r=0$ for other rays, is a vector in $\cD$.  This deformation
	is illustrated by dashed lines in the figure.  If we choose this deformation
	as $\cent\xi$, then the Scherk ends corresponding to rays $1$ and $4$ bend
	away from the parallel Scherk ends, and the embeddedness follows from
	Theorem~\ref{thm:embed2}. \qed
\end{example}

\begin{figure}[!htb]
	\includegraphics[width=.4\textwidth]{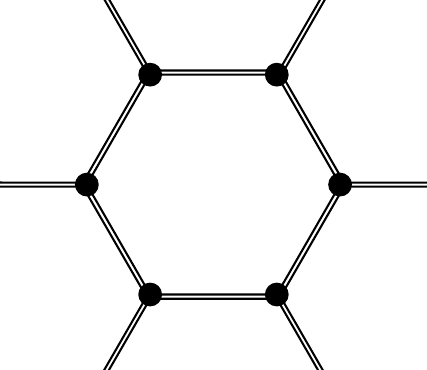}
	\caption{This graph is not rigid.\label{fig:benzene}}
\end{figure}

\begin{example}[Benzene, inconclusive]\label{ex:benzene}
	It is tempting to construct SPMSs from the Benzene-like graph in
	Figure~\ref{fig:benzene} with $12$ rays.  However, Mathematica reports that
	the space $\cD$ of balance-preserving deformations is of dimension $11 > 12 -
	2$, so the graph is not rigid.  To explicitly count the dimension of $\cD$,
	note that rotation contributes one dimension.  Up to rotation, all the
	infinitesimal deformations in $\cD$ must preserve the directions of closed
	edges; such deformations for the hexagon contribute four dimensions
	(including the scaling).  Finally, the deformations that open up parallel
	rays also preserve the balance to the first order; they contribute six
	dimensions, one for each pair.  Our main theorem is therefore inconclusive
	here.  \qed
\end{example}

\section{Construction}\label{sec:proof}

This section is dedicated to the proof of Theorem~\ref{thm:main} and the
embeddedness statements.  As we have explained, the construction of the SPMSs
will only be sketched.  The readers are referred to the first paper of this
series~\cite{saddle1} for omitted technical details.  Only the embeddedness
Theorem~\ref{thm:embed4} will receive an elaborated proof in
Section~\ref{sss:subtle} because it is our major technical concern, and the
involved argument was not detailed before.

\subsection{Weierstrass parameterization}

We construct a conformal minimal immersion using the Weierstrass
parameterization
\[
	\Sigma \ni z \mapsto \re \int^z (\Phi_1, \Phi_2, \Phi_3),
\]
where $\Sigma$ is a Riemann sphere and $\Phi_i$ are meromorphic 1-forms on
$\Sigma$ satisfying the conformality equation
\begin{equation}
	Q = \Phi_1^2 + \Phi_2^2 + \Phi_3^2 = 0. \label{eq:conformal}
\end{equation}

\subsubsection{Riemann surface} \label{ssec:opennode}

To each vertex $v \in \gV$, we associate a Riemann sphere $\hat\C_v$.  To each
half-edge $\eta \in v$, we associate a complex number $\cent p_\eta \in
\hat\C_v$, so that $\hat\C_v$ punctured at $\cent p_\eta$ provides a conformal
model for the saddle tower $\sS_v$.  Then for every $h \in \gH \setminus \gR$,
we identify $\cent p_h$ and $\cent p_{-h}$.  The resulting singular Riemann
surface with nodes is denoted $\Sigma_0$.

As $\varepsilon$ increases, we open nodes into necks in the standard way:

For each $\eta \in \gH$, let $p_\eta \in \hat\C_{v(\eta)}$ be a complex
parameter in the neighborhood of $\cent p_\eta$, and consider an adapted local
coordinate $w_h$ in a neighborhood of $p_\eta \in \hat\C_{v(\eta)}$ such that
$w_\eta(p_\eta) = 0$.  Since the graph is finite, it is possible to fix a small
number $\delta > 0$ independent of $v$ such that, for $p$ sufficiently close to
$\cent{p}$, the disks $|w_\eta|<2\delta$ are disjoint.

In the neighborhood of $0$, consider $t = (t_h)_{h \in \gH \setminus \gR}$ that
is symmetric in the sense that $t_h = t_{-h}$.  Then for every $h \in \gH
\setminus \gR$, we remove the disk
\[
	|w_h| < |t_h|/\delta,
\]
and identify the annuli
\[
	|t_h| / \delta \le |w_h| \le \delta
	\quad \text{and} \quad
	|t_{-h}| / \delta \le |w_{-h}| \le \delta
\]
by
\[
	w_h w_{-h} = t_h.
\]
This produces a Riemann surface, possibly with nodes, denoted by $\Sigma_t$,
depending on the parameters $t$ and $p$.


\subsubsection{Weierstrass data}

Let $A_h$ denote the counterclockwise circle $|w_h| = \delta$.  We need to solve the
A-period problem
\[
	\re \int_{A_\eta} (\Phi_1, \Phi_2, \Phi_3) = (0, 0, 2\pi\sigma_\eta), \quad
	\forall \eta \in \gH.
\]
For this purpose, we define $\Phi_1$, $\Phi_2$, and $\Phi_3$ as the unique
regular 1-forms on $\Sigma_t$ with simple poles at $p_r$, $r \in \gR$, possibly
also at $\infty_v$, $v\in\gV$, and the A-periods
\[
	\int_{A_\eta} (\Phi_1, \Phi_2, \Phi_3) =
	2 \pi \ii (\alpha_\eta, \beta_\eta, \gamma_\eta - \ii \sigma_\eta), \quad\forall \eta \in \gH,
\]
where $(\alpha,\beta,\gamma) \in (\cA \times \cR)^3$.  Then the A-period
problems are solved by definition.  We choose the following central value for
the parameters:
\[
	\cent\alpha_\eta = -\cos(\cent\theta_\eta),\qquad
	\cent\beta_\eta = -\sin(\cent\theta_\eta),\qquad
	\cent\gamma_\eta = 0.
\]
Then at $\varepsilon=0$ and the central values of all parameters,
$(\cent\Phi_1, \cent\Phi_2, \cent\Phi_3)$ is precisely the Weierstrass data of
the saddle tower $\sS_v$.

\subsubsection{Balance and period problems}

We want $\infty_v$ to be regular points for all $v \in \gV$.  For this purpose,
we need to solve the \emph{balance equations}
\begin{equation}
	F_v(\alpha, \beta, \gamma) =
	\sum_{\eta \in v} (\alpha_\eta, \beta_\eta, \gamma_\eta) = 0,
	\qquad \text{for all } v \in \gV. \label{eq:balance}
\end{equation}

A half-edge $r \in \gR$ corresponds to a Scherk end asymptotic to a vertical
plane. So we require that
\begin{equation} \label{eq:scherkends}
	\gamma_r \equiv 0 \quad \text{and} \quad
	|\alpha_r + \ii\beta_r| = \alpha_r^2 + \beta_r^2 \equiv 1
\end{equation}
no matter the value of other parameters.  This guarantees that (the
stereographic projection of) the Gauss map $G = -(\Phi_1 + \ii
\Phi_2)/\Phi_3$ extends holomorphically to the punctures $p_r$ with unitary
values.

Recall from~\cite{saddle1} that, for every vertex $v$, we fix an origin $O_v
\in \C_v$ bounded away from all punctures $p_h$, and a path $B_h$ from
$O_{v(h)}$ to $O_{v(-h)}$ through the neck corresponding to $h$;
see~\cite{saddle1} for the rigorous descriptions.  Then for a cycle $c=(h_1,
\cdots, h_n)$, we define $B_c$ as the concatenation $B_{h_1}*\cdots*B_{h_n}$.
For each $c \in \gF$, we need to solve the following B-period problems
\begin{equation}\label{eq:Bc}
	\re \int_{B_c} (\Phi_1, \Phi_2, \Phi_3) = 0.
\end{equation}

\subsection{Using the Implicit Function Theorem}

\subsubsection{Solving conformality problems}

At $\varepsilon=0$ and the central value of all parameters, $\cent\Phi_1$ has
$\deg(v)-2$ zeros denoted $\cent z_{v,j}$ for $1\leq j\leq \deg(v)-2$.  We may
assume that these zeros are simple and not at $\infty_v$.  When the parameters
are close to their central values, $\Phi_1$ has a simple zero $ z_{v,j}$ close
to $\cent z_{v,j}$ in $\hat\C_v$.  The same argument as in~\cite{saddle1}
proves that the conformality condition~\eqref{eq:conformal} is satisfied
if~\eqref{eq:balance} and the following equations are solved:
\begin{align}
	\int_{A_\eta} \frac{Q}{\Phi_1} &= 0, && \eta \in \gH, \label{eq:conform1} \\
	\Res{\frac{Q}{\Phi_1}}{z_{v,j}} &= 0, && 1 \le j \le \deg(v)-3,\; v \in \gV. \label{eq:conform2}
\end{align}

We make the change of parameters
\[
	\alpha_\eta + \ii\beta_\eta= - \rho_\eta \exp(\ii \theta_\eta)
\]
for $\eta \in \gH$.  The central value of $\theta_\eta$ is $\cent\theta_\eta$
given by the graph.  Recall that we write $\theta=(\theta_\eta)_{\eta \in
\gH}$.  
\begin{proposition}\label{prop:rhogammap}
	For $(t,\theta)$ in a neighborhood of $(0, \cent\theta)$, there exist unique
	values for $p_\eta$, $\rho_\eta$, and $\gamma_\eta$, depending
	real-analytically on $(t,\theta)$, such that the
	equations~\eqref{eq:conform1} and~\eqref{eq:conform2} are solved under the
	condition~\eqref{eq:scherkends}.  At $\varepsilon=0$ and the central values
	of the parameters, we have
	\[
		p_\eta = \cent p_\eta, \quad \rho_\eta = 1, \quad \gamma_\eta =0
	\]
	for $\eta \in \gH$ no matter the values of other parameters.
	Moreover, at $(t,\theta)=(0,\cent\theta)$, we have the Wirtinger derivatives
	\begin{equation}\label{eq:rhogammap}
		\frac{\partial\rho_h}{\partial t_h} = -\frac{1}{2}\Upsilon_h\Upsilon_{-h}
		\quad\text{and}\quad 
		\frac{\partial\gamma_h}{\partial t_h}=-\frac{\ii}{2}\sigma_h\Upsilon_h\Upsilon_{-h}
	\end{equation}
	for each $h \in \gH \setminus \gR$.
\end{proposition}

The proof in~\cite{saddle1} applies here almost word for word, so we omit the
proof.  From now on, we assume that $p$, $\rho$ and $\gamma$ are given by
Proposition~\ref{prop:rhogammap}.

\subsubsection{Solving horizontal balance and period problems}

For $h \in \gH \setminus \gR$, we make the change of parameters
\[
	t_h = -\exp\Big(-\ell_h\varepsilon^{-2}-\ii \sigma_h \phi_h\Big).
\]
Note that $t_h$ is a flat function in $\varepsilon$ in the sense that all
derivatives of $t_h$ in $\varepsilon$ vanish.  The central value of $\ell_h$ is
$\cent\ell_h$, the length of the segment $\varrho(e(h))$.  The central value of
$\phi_h$ is the prescribed phase function $\cent\phi_h$.  We combine $\ell$ and
$\theta$ into
\[
	x_h = \ell_h e^{\ii\theta_h},
\]
whose central value is $\cent x_h$ as given by the graph.  Recall that we write
$\vartheta=(\theta_r)_{r \in \gR}$, and $\chi=(x, \vartheta)$.  For $\chi$ in
an neighborhood of $\cent\chi$, we change to the variable
\[
	\chi = \cent\chi + \varepsilon^2\xi + \zeta
\]
where $\xi \in \cD$ and $\zeta \in \cD^\perp$.

\begin{proposition}\label{prop:thetaell}
	Assume that the graph $\Gamma$ is balanced and rigid.  For
	$(\varepsilon,\xi,\phi)$ in a neighborhood of $(0, \cent\xi, \cent\phi)$,
	there exist unique values for $\zeta$, depending smoothly on $(\varepsilon,
	\xi, \phi)$, such that the horizontal components of the balance
	equations~\eqref{eq:balance} and the B-period equations \eqref{eq:Bc} are
	solved.  Moreover, $\zeta$ is an even function of $\varepsilon$ and, at
	$\varepsilon=0$, we have $\zeta(0,\xi,\phi) = 0$ and
	$\dfrac{1}{2}\dfrac{\partial^2\zeta}{d\varepsilon^2}(0,\xi,\phi) = \dot\zeta$
	is the unique solution to~\eqref{eq:system} in $\cD^\perp$ no matter the
	values of $\xi$ and $\phi$.
\end{proposition}

The proof in~\cite{saddle1} applies here with some modification. We sketch a
proof here because some computations will be useful later for the embeddedness
proofs.

\begin{proof}[Sketched proof]
	Define for $(\varepsilon,\xi,\phi)$ in a neighborhood of
	$(0,\cent\xi,\cent\phi)$ and $h \in \gH \setminus \gR$
	\[
		\cPhor_h(\varepsilon,\xi,\phi)
		= \varepsilon^2\Big(\re \int_{B_h} \Phi_1 + \ii \re \int_{B_h} \Phi_2\Big),
	\]
	for $c \in \gF$
	\[
		\cPhor_c(\varepsilon,\xi,\phi)
		= \varepsilon^2\Big(\re \int_{B_c} \Phi_1 + \ii \re \int_{B_c} \Phi_2\Big)
		= \sum_{h \in c} \cPhor_h(\varepsilon,\xi,\phi),
	\]
	and for $v \in \gV$
	\begin{equation}\label{eq:hforce}
		\cFhor_v(\varepsilon,\xi,\phi) = - \div_v(\alpha + \ii\beta)
		= \sum_{h \in v} \rho_h e^{\ii\theta_h}.
	\end{equation}

	In~\cite{saddle1} we have computed that
	\begin{equation}\label{eq:hperiod}
		\cPhor_h(\varepsilon,\xi,\phi) =
		\rho_h x_h + \varepsilon^2 \lambda_h(\varepsilon,\xi,\phi)
	\end{equation}
	where $\lambda_h$ is analytic in $t_h$ and $\lambda_h(0,\xi,\phi) = \mu_h -
	\mu_{-h} =: \mu^a_h$ no matter the values of $\xi$ and $\phi$.  So
	$\lambda_h$ is flat in $\varepsilon$.  We have also computed that
	\begin{equation}\label{eq:systemanalnew}
		\begin{cases}
			\cPhor(\varepsilon,\xi,\phi) = \Phor(\chi) + \varepsilon^2 \Phor(\mu^a) + \text{flat terms},\\
			\cFhor(\varepsilon,\xi,\phi) = \Fhor(\chi) + \text{flat terms}.
		\end{cases}
	\end{equation}

	Write $\cFhor = (\cFhor_v)_{v \in \gV}$ and $\cPhor = (\cPhor_c)_{c \in
	\gF}$.  We want to solve
	\begin{equation}\label{eq:horsystem}
		(\cFhor, \cPhor)(\varepsilon,\xi,\phi) = 0.
	\end{equation}
	If $\Gamma$ is balanced, the system is solved at $\varepsilon=0$ with
	$\zeta=0$ no matter the values of $\xi$ and $\phi$.  If $\Gamma$ is rigid,
	then by the Implicit Function Theorem, the system has a unique solution
	$\zeta(\varepsilon,\xi,\phi)$, depending smoothly on $(\varepsilon,\xi,\phi)$
	in a neighborhood of $(0, \cent\xi, \cent\phi)$.

	The system~\eqref{eq:horsystem} is even in $\varepsilon$, so must be the
	solution $\zeta$.  Taking the second derivative of~\eqref{eq:horsystem} with
	respect to $\varepsilon$ at $\varepsilon = 0$ gives the linear system
	\[
		\begin{dcases}
			2\Phor_c(\mu^a) + \Phor_c\Big(\frac{\partial^2\zeta}{\partial\varepsilon^2}\Big) = 0, &c \in \gF,\\
			D\Fhor_v(\cent\chi) \cdot \Big(\frac{\partial^2\zeta}{\partial\varepsilon^2}\Big) = 0, & v \in \gV,
		\end{dcases}
	\]
	which proves that
	$\dfrac{1}{2}\dfrac{\partial^2\zeta}{\partial\varepsilon^2}$ must be the
	unique solution $\dot\zeta \in \cD^\perp$ to the linear
	system~\eqref{eq:system} in $\cD^\perp$.
\end{proof}

\subsubsection{Solving vertical balance and period problems}

\begin{proposition}\label{prop:phi}
	Assume that the phase function $\cent\phi$ is balanced and rigid with respect
	to the prescribed deformation $\cent\xi \in \cD$.  For $(\varepsilon,\xi)$ in
	a neighborhood of $(0,\cent\xi)$, there exist unique values for $(\phi_h)_{h
	\in \gH\setminus\gR}$, depending smoothly on $\varepsilon$ and $\xi$, such
	that $\phi_h(0,\cent\xi) = \cent\phi_h$, and the vertical components of the
	balance equations~\eqref{eq:balance} and the B-period equations~\eqref{eq:Bc}
	are solved.
\end{proposition}

The proof in~\cite{saddle1} applies here with only slight modification, but we
still sketch a proof for completeness.

\begin{proof}[Sketched proof]
	Define for $(\varepsilon,\xi)$ in a neighborhood of $(0,\cent\xi)$ and $c \in
	\gF$
	\[
		\cPver_c (\varepsilon,\xi,\phi)= \re \int_{B_c} \Phi_3,
	\]
	and for $\varepsilon > 0$ and $b \in \gB_m^*$,
	\[
		\cFver_b(\varepsilon,\phi) := -\exp(\cent\ell_b\varepsilon^{-2}) \div_b(\gamma(\varepsilon,\xi,\phi)).
	\]

	By the same computation as in~\cite{saddle1}, we have
	\[
		\cPver_c(0,\xi,\phi) = \Pver_c(\phi) \pmod{2\pi}
	\]
	no matter the value of $\xi$, and that $\cFver_b$ extends smoothly at
	$\varepsilon=0$ to
	\[
		\cFver_b(0,\xi,\phi) = \sum_{h\in m(b)}\Upsilon_h\Upsilon_{-h}\sin(\phi_h)\exp\left(-\re({\dot x}_h\exp(-\ii\cent\theta_h))\right),
	\]
	where $\xi+\dot\zeta = (\dot x, \dot \vartheta)$.  In particular,
	\[
		\cFver_b(0,\cent\xi,\phi) = \Fver_b(\phi).
	\]

	Write $\cFver = (\cFver_b)_{b \in \gB_m^*}$ and
	$\cPver = (\cPver_c)_{c \in \gF}$.  We want to solve
	\[
		(\cFver, \cPver)(\varepsilon,\xi,\phi) = 0.
	\]
	If the phase function $\cent\phi$ is balanced with respect to $\cent\xi$, the
	system is solved at $(\varepsilon,\xi)=(0,\cent\xi)$ by $\phi=\cent\phi$.  If
	$\cent\phi$ is rigid, then by the Implicit Function Theorem, the system has a
	unique solution $\phi(\varepsilon,\xi)$, depending smoothly on
	$(\varepsilon,\xi)$ in a neighborhood of $(0,\cent\xi)$, such that
	$\phi(0,\cent\xi) = \cent\phi$.
\end{proof}

The same argument as in~\cite{saddle1} shows that the immersion is regular.
This finishes the proof of Theorem~\ref{thm:main}.

\subsection{Symmetric SPMSs}

A trivial phase function is trivially balanced, but not necessarily rigid.  So
Theorem~\ref{thm:symmetry} is not contained in our main Theorem~\ref{thm:main}.
But its proof is very similar, only much easier, hence we only give a brief
sketch here.  The readers are referred to~\cite{younes2009}
and~\cite{traizet2001} for technical details.

The reflection in the horizontal symmetry plane correspond to an involution
$\rho$ of $\Sigma_t$ that restricts to $\rho(z) = \overline z$ on $\hat\C_v$
for every vertex $v$.  We restrict to the parameters $p_h$ to real values, and
set $\gamma \equiv 0$ so that the vertical balance problem is trivially solved.
This ensures that $\rho^*\Phi_{1,2} = \overline\Phi_{1,2}$ and $\rho^*\Phi_3 =
-\overline\Phi_3$, so the surface carries the desired symmetry;
see~\cite{younes2009}.

Since $\phi \equiv \cent\phi$ is trivial, we have $t_h \in \R$, negative if
$\phi_h = 0$, positive if $\phi_h = \pi$.  For each $c \in \gF$, we choose the
integral path $B_c$ as the concatenation of $B'_h$ where
\begin{itemize}
	\item if $\phi_h = 0$, $B'_h$ consists of the real segment $w_h = -t_h/\delta >
		0$ to $w_{\varsigma(h)} = -\delta < 0$.

	\item if $\phi_h = \pi$, $B'_h$ consists of an clockwise half-circle around
		$p_h$ from $w_h = -t_h/\delta < 0$ to $w_h = t_h/\delta > 0$, followed by
		the real segment to $w_{\varsigma(h)} = -\delta$.
\end{itemize}
This careful choice of path  makes it convenient to compute that $\re\int_{B_c}
\Phi_3 = 0 \pmod{2\pi}$; see~\cite{younes2009}.  So the vertical period problem
is automatically solved.

We then use the Implicit Function Theorem to solve the conformality problem and
horizontal balance and period problems, as we did above for the general case.
The result is a continuous family of symmetric SPMSs depending on $|\gR|-1$
parameters.  One of them is $\varepsilon$.  The other parameters correspond to
the deformations $\chi \in \cM$; cf.~\cite{traizet2001}.  This finished the
construction of symmetric SPMSs.

\subsection{Embeddedness}\label{ssec:embed}

We now prove the scenarios where the surfaces $M_{\varepsilon,\xi}$ in
Theorem~\ref{thm:main} are embedded, at least for some $(\varepsilon,\xi)$ in a
neighborhood of $(0,\cent\xi)$.  In fact, the same proof as in~\cite{saddle1}
proves the embeddedness except for the Scherk ends.  It also proves the
embeddedness of each Scherk end.  The only problem is that the Scherk ends
might intersect each other.  We then omit this part of the proof, and focus on
the bendings of the Scherk ends.

\subsubsection{The safe case}

If the graph has no parallel rays, then the Scherk ends would not intersect for
$(\varepsilon,\xi)$ sufficiently close to $(0,\cent\xi)$.  This proves the
scenario in Theorem~\ref{thm:embed1}, and has been considered
in~\cite{traizet2001}.

\subsubsection{Lowest-order bending}

Otherwise, we must analyse how the Scherk ends bend.  For this purpose, let us
consider a 1-parameter family $M_{\varepsilon,\xi(\varepsilon)}$ where
$\xi(\varepsilon) \in \cD$ is a fixed analytic function such that $\xi(0) =
\cent\xi$.  All other parameters have been solved as a smooth function of
$\varepsilon$.  Write
\[
	\zeta(\varepsilon) = \widetilde\zeta(\varepsilon) +
	\overline\zeta(\varepsilon).
\]
Here and in the remaining of the paper, for any smooth function $f(x)$, we use
$\widetilde f(x)$ to denote the analytic function given by the Taylor series of
$f$ in $x$, and use $\overline f(x)$ to denote the non-analytic remainder.

By Proposition~\ref{prop:thetaell}, the first term of
$\widetilde\zeta(\varepsilon)$ is $\varepsilon^2 \dot\zeta$.  Hence the
lowest-order deformation of the graph is $\cent{\dot\chi} =
\cent\xi+\dot\zeta$.  If the parallel Scherk ends bend away from each other
under this deformation, the Scherk ends would not intersect for
$(\varepsilon,\xi)$ sufficiently close to $(0,\cent\xi)$.  This proves the
scenario in Theorem~\ref{thm:embed2}.

\subsubsection{Analytic bending}

If the lowest-order deformation does not help, we may look into higher order
deformations, and eventually use the entire Taylor series.  In view
of~\eqref{eq:systemanalnew}, $\widetilde\zeta(\varepsilon) \in \cD^\perp$
solves the non-linear system
\[
	\begin{cases}
		\widetilde{\cPhor_c} = \Phor_c(\widetilde\zeta) + \varepsilon^2 \Phor_c(\mu^a_h) = 0, & c \in \gF,\\
		\widetilde{\cFhor_v} = \Fhor_v(\widetilde\chi) = 0, & v \in \gV,
	\end{cases}
\]
where $\widetilde\chi = \cent\chi + \varepsilon^2 \xi(\varepsilon) +
\widetilde\zeta(\varepsilon)$.  This is exactly the
system~\eqref{eq:systemanal}.  Moreover, since the graph is balanced and rigid,
there is a unique solution $\widetilde\zeta(\varepsilon)$ for $\varepsilon$
sufficiently small.  If the parallel Scherk ends bend away from each other
under the deformation $\widetilde\chi(\varepsilon)$ for sufficiently small
$\varepsilon$, then the Scherk ends would not intersect for sufficiently small
$\varepsilon$.  This proves the scenario described in Theorem~\ref{thm:embed3}.

\subsubsection{Non-analytic bending}\label{sss:subtle}

The remaining of the paper is dedicated to the proof of
Theorem~\ref{thm:embed4}.  This is the only part of the paper that contains
detailed arguments, because the technical details here were never written down
before.

\medskip

Interestingly, if the graph appears as a simple line arrangement, the analytic
part $\widetilde\chi(\varepsilon)$ does not resolve parallel Scherk ends no
matter how many terms are used.  This follows from the following lemma.

\begin{lemma}\label{lem:lines}
	Write $\widetilde\chi(\varepsilon) = (\widetilde x(\varepsilon),
	\widetilde\vartheta(\varepsilon) )$.  If the graph appears as a simple line
	arrangement, then we have $\widetilde\theta_r = \widetilde\theta_{r'} + \pi
	\pmod{2\pi}$ whenever the rays $r$ and $r'$ belong to the same line.
\end{lemma}

\begin{proof}
	Explicitly, we have the Taylor series
	\[
		\widetilde{\cFhor_v} = \Fhor_v(\widetilde\chi) = - \sum_{\eta \in
		v} e^{\ii \widetilde\theta_\eta}.
	\]
	Since the graph appears as a simple line arrangement, every vertex is of
	degree four.  Consequently, $\widetilde{\cFhor}=0$ if and only if
	\[
		\widetilde\theta_{\rotate^2(\eta)} = \widetilde\theta_\eta + \pi \pmod{2\pi}.
	\]
	So, if two Scherk ends correspond to rays from the same line, their
	directions $\widetilde\theta$ must remain opposite.
\end{proof}

Now assume that the line arrangement contains a pair of parallel rays.  If they
bend away from each other under the deformation $\widetilde\chi$, then by the
Lemma, the other rays on the same lines must bend towards each other, creating
an unwanted intersection.  The only way to avoid this is to let the parallel
rays remain parallel under the deformations $\widetilde\chi$, i.e.
$\widetilde\theta_r(\varepsilon) = \widetilde\theta_{r+1}(\varepsilon)$ for
sufficiently small $\varepsilon$ whenever $\cent\theta_r = \cent\theta_{r+1}$.
This is assumed in Theorem~\ref{thm:embed4}.

\begin{remark}
	There are certainly other situations that $\widetilde\chi$ can not resolve
	parallel Scherk ends; see Example~\ref{ex:paralleltree}.  But we have no plan
	to classify all such situations.
\end{remark}

Then we proceed to investigate the non-analytic part
$\overline\zeta(\varepsilon)$.  Let $\cent \ell_{\min}$ be the length of the
shortest edges in the graph, and write $\tau(\varepsilon) = \exp(-\cent
\ell_{\min}/\varepsilon^2)$.  Recall that $m(\gH) = \{ h \in \gH \setminus \gR
\colon \cent\ell_h = \cent\ell_{\min} \}$.

\begin{proposition}
	\[
		\overline\zeta \sim \widehat\zeta \tau
	\]
	where $\widehat\zeta$ is the unique solution in $\cD^\perp$
	to~\eqref{eq:systemflat}.
\end{proposition}

\begin{proof}
	We first prove that $\overline\zeta \in O(\tau)$.  Assume instead that
	$\overline\zeta \sim \kappa f(\varepsilon)$ for some $\kappa \in \cD^\perp$
	and $f(\varepsilon) \in \omega(\tau(\varepsilon))$ ($f$ dominates $\tau$) as
	$\varepsilon \to 0$.  In~\cite{saddle1}, we have computed
	from~\eqref{eq:rhogammap} that
	\[
		\overline\gamma_h \sim -K_h \sin\cent\phi_h e^{-\cent\ell_h/\varepsilon^2}.
	\]
	A similar computation yields
	\[
		\overline\rho_h \sim K_h \cos\cent\phi_h e^{-\cent\ell_h/\varepsilon^2}.
	\]
	Moreover, we have seen that $\lambda_h(\varepsilon)$ is analytic in $t_h$, so
	$\overline\lambda_h \in O(e^{-\cent\ell_h/\varepsilon^2})$.  As $\cent\ell_h \ge
	\cent\ell_{\min}$, the functions $\overline\rho$ and $\overline\lambda$ are
	all in $O(\tau)$, hence all dominated by $f$ as $\varepsilon \to 0$.  Then a
	routine computation from~\eqref{eq:hforce} and~\eqref{eq:hperiod} yields, as
	$\varepsilon \to 0$, that
	\[
		\begin{cases}
			f(\varepsilon)^{-1}\overline{\cPhor}(\varepsilon) \to
			\Phor_c(\kappa) = 0, & c \in \gF,\\
			f(\varepsilon)^{-1}\overline{\cFhor}(\varepsilon) \to
			D\Fhor_v(\cent\chi) \cdot \kappa = 0, & v \in \gV,
		\end{cases}
	\]
	hence $\kappa = 0$.

	Then we compute the coefficient $\widehat\zeta \in \cD^\perp$ of $\tau$.
	Note that $\overline\rho_h \in o(\tau)$ as $\varepsilon \to 0$ if and only if
	$\cent\ell_h > \cent \ell_{\min}$, i.e.\ $h \not\in m(\gH)$.  Using this
	fact, a routine computation from~\eqref{eq:hforce} and~\eqref{eq:hperiod}
	yields, as $\varepsilon \to 0$, that
	\[
		\begin{dcases}
			\tau(\varepsilon)^{-1} \overline{\cPhor}(\varepsilon) \to
			\Phor_c(\widehat\zeta)
			+ \sum_{h \in c \cap m(\gH)} x_h K_h \cos \cent\phi_h
			= 0, & c \in \gF,\\
			\tau(\varepsilon)^{-1} \overline{\cFhor}(\varepsilon) \to
			D\Fhor_v(\cent\chi) \cdot \widehat\zeta
			+ \sum_{h \in v \cap m(\gH)} u_h K_h \cos \cent\phi_h
			= 0, & v \in \gV,
		\end{dcases}
	\]
	so $\widehat\zeta$ must be the unique solution in $\cD^\perp$
	to~\eqref{eq:systemflat}.

	Finally, we prove that $\widehat\zeta \ne 0$.  This follows from the
	existence of a vertex $v$ such that the summation
	\begin{equation}\label{eq:vsum}
		\sum_{h \in v \cap m(\gH)} u_h K_h \cos \phi_h
	\end{equation}
	is not zero.

	Assume the opposite, i.e.\ that the summation vanishes for every vertex.
	Define
	\[
		V_m = \{ v \in V \colon \exists h \in v\, \text{such that}\, \cent\ell_h =
		\cent\ell_{\min}\,\text{and}\,K_h \cos\cent\phi_h \ne 0\}.
	\]
	Clearly, if $V_m$ is not empty, it must contain at least two vertices.  Note
	that $v \cap m(\gH) = m(v)$ if $v \in V_m$.

	Take the convex hull of $\{\varrho(v) \colon v \in V_m\}$, and consider an
	arbitrary vertex $v$ at a corner of the convex hull.  Since every vertex is
	of degree four, $m(v)$ consists of either a single half-edge, or two
	half-edges whose corresponding unit vectors are linearly independent.  In
	either case, the summation~\eqref{eq:vsum} vanishes if and only if $K_h
	\cos\cent\phi_h = 0$ for all $h \in m(v)$.  This contradicts our assumption
	that $v \in V_m$.

	So $V_m$ is empty, meaning that $K_h \cos\cent\phi_h = 0$ for all $h \in
	m(\gH)$.  Recall from~\cite{saddle1} that $K_h > 0$, so $\cos\cent\phi_h = 0$
	for all $h \in m(\gH)$.  Recall that the vertical force along $h$ is $K_h
	\sin\cent\phi_h$.  So $K_h \cos\cent\phi_h$ is nothing but the derivative of
	the vertical force with respect to $\phi_h$ at $\cent\phi_h$.  Hence for a
	cut $b \in \gB_m$ such that $\cent\ell_b = \cent\ell_{\min}$, the derivative
	of $\Fver_b$ with respect to $\phi$ at $\cent\phi$ is $0$.  This contradicts
	the assumption of Theorem~\ref{thm:main} that the phase function is rigid.

	This finishes the proof that the summation \eqref{eq:vsum} must be non-zero
	for some vertex.  We then conclude that $\widehat\zeta \ne 0$.
\end{proof}

If the graph appears as a simple line arrangement, and the deformation
$\widetilde\chi$ does not create any self-intersection, then parallel Scherk
ends are resolved if they bend away from each other under the deformation
$\widehat\zeta$.  This finises the proof for Theorem~\ref{thm:embed4}.

\appendix

\section{Rigidity for simple graphs}

We now prove Proposition~\ref{prop:rigid}.  For this purpose, let us
restate~\cite[Lemma 5.1]{perez2007} in the following form.

\begin{lemma}
	Under the conditions of Proposition~\ref{prop:rigid}, let $v$ be an arbitrary
	vertex and $L$ be a straight line though $\rho(v)$ that does not contain the
	geometric representation of any edge adjacent to $v$, then each side of $L$
	contains the geometric representations of at least two edges adjacent to $v$.
\end{lemma}

This lemma allows us to adapt the argument in~\cite{traizet2001} that proves
the rigidity when the graph appears as a line arrangement.

\begin{proof}[Sketched proof of Proposition~\ref{prop:rigid}]
	Up to a rotation, we may assume that the vectors $x_h \in \C$, $h \in \gH$,
	all have non-zero real-parts.  As a consequence, the geometric
	representations of the vertices all have distinct real parts.  We then order
	the vertices by the real parts of $\rho(v)$, and order the faces by their
	left-most vertices.

	For each vertex $v \in \gV$, we choose two half-edges in $v$ whose
	directions, say $\cent\theta_{v,1}$ and $\cent\theta_{v,2}$, point to the
	left side; their existence follows from the lemma above.  Then $D\Fhor$
	restricted to the variables $\dot\theta_{v,1}$ and $\dot\theta_{v,2}$ is a
	real square matrix with $|\gV| \times |\gV|$ blocks of size $2 \times 2$.
	All blocks above the diagonal are $0$.  To see this, note that if $v' < v$,
	then $\Fhor_v$ is independent of the chosen half-edges in $v'$.  The diagonal
	blocks are invertible because the unit vectors $e^{\ii\cent\theta_{v,1}}$ and
	$e^{\ii\cent\theta_{v,2}}$ are linearly independent; see~\cite{traizet2001}.
	This proves that $D\Fhor$ is surjective.

	For each face $c \in \gF$, we choose two half-edges in $c$ whose edges are
	adjacent to the left-most vertex of $c$. Let $\cent\ell_{c,1}$ and
	$\cent\ell_{c,2}$ be their lengths.  Then $D\Phor$ restricted to the
	variables $\dot\ell_{c,1}$ and $\dot\ell_{c,2}$ is a real square matrix with
	$|\gF| \times |\gF|$ blocks of size $2 \times 2$.  All blocks above the
	diagonal are $0$.  To see this, note that if $c' < c$, then $\Phor_c$ is
	independent of chosen half-edges in $c'$.  The diagonal blocks are invertible
	because the unit vectors $e^{\ii\cent\theta_{c,1}}$ and
	$e^{\ii\cent\theta_{c,2}}$ are linearly independent; see~\cite{traizet2001}.
	This proves that $D\Phor$ is surjective.

	Finally, $(D\Fhor, D\Phor)$ restricted to the variables $\dot x_{v,1}$, $\dot
	x_{v,2}$, $\dot x_{c,1}$, $\dot x_{c,2}$, where $v \in \gV$ and $c \in \gF$,
	is a real square matrix. It can be partitioned into four blocks. The two
	diagonal blocks, of size $2|\gV| \times 2|\gV|$ and $2|\gF| \times 2|\gF|$
	respectively, are invertible by the argument above.  Clearly, $\Fhor$ does
	not depend on $\ell$, hence one off-diagonal block is zero.  The matrix is
	therefore invertible; see~\cite{traizet2001}.  This proves the rigidity of
	the graph.
\end{proof}

\begin{remark}
	In the presence of parallel edges, the proof fails because the vectors of
	chosen half-edges might be linearly dependent, hence the diagonal blocks
	might be singular.
\end{remark}

\bibliography{References}
\bibliographystyle{alpha}

\end{document}